\documentclass{article}
\usepackage{amsopn,amstext,amsbsy,amsmath,amscd,amsthm,latexsym,amsfonts}
\usepackage{geometry} 
\usepackage{graphicx}
\usepackage[utf8]{inputenc} 
\usepackage{epstopdf}
\usepackage{amssymb}
\usepackage{mathtools}
\usepackage{graphics}
\usepackage{pgfplots}
\usepackage[hang,centerlast]{subfigure}
\usepackage{float}
\usepackage[titletoc,title]{appendix}
\usepackage[swedish,german,english]{babel}
\pgfplotsset{compat=1.12}
\usepackage[all]{xy} 
\allowdisplaybreaks[1]
\theoremstyle{plain}
\newtheorem{theorem}                 {\bf Theorem}

\newtheorem{lemma}        [theorem]  {\bf Lemma}

\title{Three-dimensional internal gravity-capillary waves in finite depth}
\author{Dag Nilsson}
\theoremstyle{definition}
\newtheorem{remark}       [theorem] {\bf Remark}

\newcommand\norm[1]{\left\lVert#1\right\rVert}

\renewcommand{\i}{{\mathrm{i}}}

\newcommand{\abs}[1]{\lvert#1\rvert}
\begin{document}
\maketitle
\begin{abstract}
We consider three-dimensional inviscid irrotational flow in a two layer fluid under the effects of gravity and surface tension, where the upper fluid is bounded above by a rigid lid and the lower fluid is bounded below by a flat bottom. We use a spatial dynamics approach and formulate the steady Euler equations as an infinite-dimensional Hamiltonian system, where an unbounded spatial direction $x$ is considered as a time-like coordinate. In addition we consider wave motions that are periodic in another direction $z$. 
By analyzing the dispersion relation we detect several bifurcation scenarios, two of which we study further: a type of $00(\mathrm{i}s)(\mathrm{i}\kappa_0)$ resonance and a Hamiltonian-Hopf bifurcation. The bifurcations are investigated by performing a center-manifold reduction, which yields a finite-dimensional Hamiltonian system. For this finite-dimensional system we establish the existence of periodic and homoclinic orbits, which correspond to, respectively, doubly periodic travelling waves and oblique travelling waves with a dark or bright solitary wave profile in the $x$-direction. The former are obtained using a variational Lyapunov-Schmidt reduction and the latter by first applying a normal form transformation and then studying the resulting canonical system of equations.
  
\end{abstract}
\section{Introduction}
\subsection{Internal waves}
Internal waves are waves which propagate along the interface of two immiscible fluids of different density. In this paper we study three-dimensional internal waves under the influence of gravity and interfacial tension. The flow is assumed to be inviscid and irrotational and the density of each layer is assumed to be constant. In addition we assume that the upper fluid is bounded above by a rigid horizontal lid and the lower fluid is bounded below by a rigid horizontal bottom. The two fluids are separated by an interface $\eta$, which is a function of $X,Z$, in the domain $\{(X,Y,Z)\in \mathbb{R}^3\ : \ -h_2\leq Y\leq h_1\}$, where $h_1,h_2$ are positive real numbers. Let $\rho_1,\rho_2$ be the densities of the upper and lower fluid respectively, where $\rho_1<\rho_2$, and let $\phi_1,\phi_2$ be the velocity potentials of the upper and lower fluid respectively. We consider waves which travel with constant speed $c$ in the positive $X$-direction. The governing equations can then be written as 
\begin{align}
\Delta \phi_1&=0,\quad \text{for } \eta<Y<h_1,\label{goveq1}\\
\Delta \phi_2&=0,\quad \text{for }-h_2<Y<\eta,
\end{align}
with boundary conditions
\begin{align}
\phi_{1Y}&=0&& \text{on } y=h_1,\\
\phi_{2Y}&=0&&\text{on } y=-h_2,\\
\phi_{1Y}&=-c\eta_X+\eta_X\phi_{1X}+\eta_Z\phi_{1Z}&&\text{on } Y=\eta,\\
\phi_{2Y}&=-c\eta_X+\eta_X\phi_{2X}+\eta_Z\phi_{2Z}&&\text{on } Y=\eta,\\
\rho_2\bigg(&-c\phi_{2X}+\frac{1}{2}\abs{\nabla\phi_2}^2+g\eta\bigg)-\rho_1\left(-c\phi_{1X}+\frac{1}{2}\abs{\nabla\phi_1}^2+g\eta\right)\nonumber\\
&=\sigma\left(\frac{\eta_X}{\sqrt{1+\eta_X^2+\eta_Z^2}}\right)_X+\sigma\left(\frac{\eta_Z}{\sqrt{1+\eta_X^2+\eta_Z^2}}\right)_Z&&\text{on }\label{goveq8} Y=\eta,
\end{align}
\begin{figure}\centering
\includegraphics{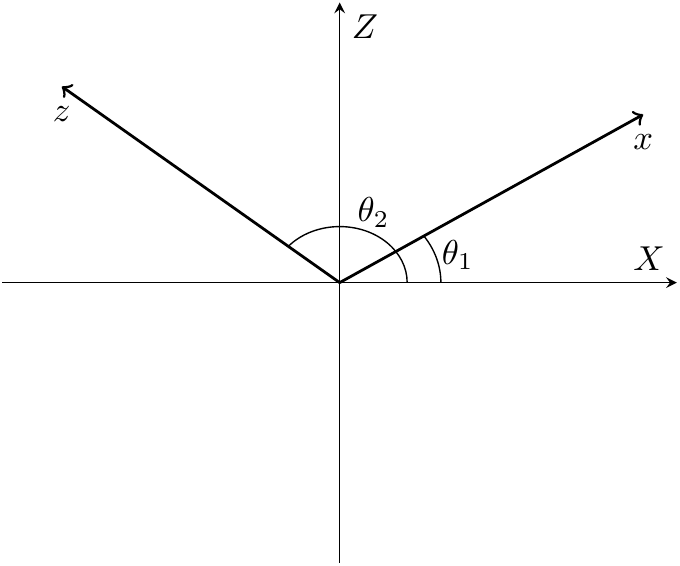}
\caption{We are considering travelling waves with a bounded profile in the $x$ direction that are periodic in $z$, with $X$ being the direction of propagation.}
\label{angles}
\end{figure}\noindent
where $\sigma$ is the coefficient of interfacial tension and $g$ is the gravitational constant. In addition we will consider waves which have a bounded profile in some direction $x$, and are periodic in some other direction $z$. Let $\theta_1$ be the angle between the $x$-axis and the $X$-axis and let $\theta_2$ be the angle between the $z$-axis and the $X$-axis (see Figure \ref{angles}), so that
\begin{equation}\label{newcoordinates}
x=\cos(\theta_1)X+\sin(\theta_1)Z,\quad z=\cos(\theta_2)X+\sin(\theta_2)Z.
\end{equation}
Solutions of \eqref{goveq1}--\eqref{goveq8} which depend upon $x,Y,z$ and are periodic in $z$, are called oblique travelling waves. Oblique travelling waves for which there exist angles $\theta_1$, $\theta_2$ such the waves are independent of $z$, are called oblique line waves. In order to find oblique travelling wave solutions we will use the method of spatial dynamics. The idea, which is due to Kirchgässner \cite{K82}, is to formulate a time-independent problem as an evolution equation in which a spatial coordinate plays the role of time. In our case we will use $x$ as time and obtain the evolution equation
\begin{equation}\label{modelhameq}
u_x=Ku+\tilde{\mathcal{F}}(u),
\end{equation}
where $u$ belongs to some Banach space, $K$ is a linear operator and $\tilde{\mathcal{F}}(u)=\mathcal{O}(\norm{u}^2)$. This is an ill-posed problem but it is possible to obtain bounded solutions by applying the center-manifold theorem. This is a result which can be used to obtain a finite-dimensional system of equations on a center manifold, which is locally equivalent to the original equation \eqref{modelhameq}. 
The idea of using spatial dynamics to study oblique travelling is due to Groves and Haragus \cite{Groves2001} and in the present work we rely heavily on the methods developed by them.
\subsection{Previous work on three-dimensional surface waves}
We mention here some relevant results concerning three-dimensional surface waves that are periodic in at least one distinguished direction.
Groves and Mielke \cite{Groves2007} considered waves that are periodic in the transverse direction $Z$, with a bounded profile in the direction of propagation $X$. This correspond to choosing $\theta_1=0$, $\theta_2=\pm\pi/2$ in our setting. More specifically they construct waves with a periodic, quasiperiodic or generalized solitary wave profile in the $X$-direction. Waves which are periodic in the direction of propagation $X$ and have a bounded profile in the transverse direction $Z$ (that is, choosing $\theta_1=\pm \pi/2$, $\theta_2=0$), were studied in \cite{Groves2001,HK2001}. The authors found waves with a periodic, quasiperiodic or generalized solitary wave profile in the $Z$-direction. The general case, that is when arbitrary angles $\theta_1,\theta_2\in (-\pi,\pi)$ are allowed, were considered in \cite{Groves2003}. Here waves which have a bounded profile in some general direction $x$ and are periodic in some other direction $z$, were found. All of the above mentioned results were obtained by applying the method of spatial dynamics, using a formulation as in \eqref{modelhameq}. In \cite{Groves2001,Groves2007,HK2001} the spectrum of $K$ depends upon the parameters $\alpha=gh/c^2$, $\beta=T/hc^2$ and $\nu=2\pi/P$, where $c$ is the wave speed, $\sigma$ is the coefficient of surface tension, $h$ is the water depth and $P$ is the period either in the direction of propagation or in the transverse direction. The parameters $\alpha$ and $\beta$ emerge when the governing equations are nondimensionalized, and $\nu$ appears when the period is normalized to $2\pi$. In \cite{Groves2003} the spectrum also depends upon the angles $\theta_1$ and $\theta_2$. These extra parameters allow for a plethora of different bifurcation scenarios. In fact, as was observed in \cite{Groves2003}, essentially all possible bifurcation scenarios known in Hamiltonian systems theory can be obtained by varying the different parameters. 
We also mention some results on doubly periodic waves obtained using other methods than spatial dynamics. Reeder and Shinbrot \cite{Reeder1981} proved the existence of doubly periodic waves with a diamond pattern, that is with $\theta_2=-\theta_1$. This is done by solving an associated linear problem and then, using the solutions of the linear problem, constructing a sequence whose limit is a solution of the full nonlinear problem. Iooss and Plotnikov \cite {IP3d} considered the same problem as in \cite {Reeder1981} but in the absence of surface tension. The absence of surface tension gives rise to a small divisor problem and the authors use Nash-Moser methods to prove the existence of doubly periodic waves with a diamond pattern.
In \cite{Craig2000} the authors proved the existence of doubly periodic waves with arbitrary angles $\theta_1$, $\theta_2$, using a variational approach and more specifically, a variational Lyapunov-Schmidt reduction.
\subsection{Previous work on three-dimensional internal waves}
Three-dimensional travelling internal waves are not as well-studied as their surface wave counterparts. In particular, to the author's knowledge there are no rigorous existence results for such waves. There are however several numerical results concerning such waves, see for example \cite{Parau2007,Parau2007a}. There is also a recent work \cite{Akers2017} where the authors study overturning waves propagating on the interface between two fluids. We also mention \cite{KIM2006} where the authors studied an extension of the Benjamin equation, which can be used to model internal waves, and were able to show that it posses fully localized solitary wave solutions.

\subsection{Outline of paper}
In section \ref{spatialdyn} the parameters $\alpha=gh_1(1-\rho)/c^2,\ \beta=\sigma/(h_1\rho_2c^2),\ \rho=\rho_1/\rho_2,\ h=h_2/h_1$ emerge from the nondimensionalization of the governing equations \eqref{goveq1}--\eqref{goveq8}. We then obtain a Hamiltonian formulation of the problem by first identifying solutions of the governing equations as critical points of a certain functional. This functional is found from Luke's variational principle and can be identified as an action integral, from which a Hamiltonian is obtained by performing a Legendre transform. 
The boundary conditions associated with the corresponding Hamiltonian system are nonlinear, whereas the center-manifold theorem applies to equations on linear spaces. It is therefore necessary to perform a change of variables so that we get a Hamiltonian system with linear boundary conditions. This is done in section \ref{changevar}. 

The dimension of the center manifold is equal to the number of imaginary eigenvalues of $K$, counted with multiplicity. Due to this we carry out an investigation of the spectrum of $K$ in section \ref{spectrum}. Since we assume periodicity in $z$ we expand in Fourier series and consider the eigenvalue equation for each Fourier mode $k$. We find that an imaginary number $\i s$ is a mode $k$ eigenvalue if and only if the dispersion relation is satisfied:
\begin{equation}\label{disprelintro}
\frac{\rho\left(\nu k\cos(\theta_2)+s\cos(\theta_1)\right)^2}{\tanh(\tilde{\gamma}_k)}+\frac{\left(\nu k\cos(\theta_2)+s\cos(\theta_1)\right)^2}{\tanh(h\tilde{\gamma}_k)}=(\alpha+\beta\tilde{\gamma}_k^2)\tilde{\gamma}_k,
\end{equation}
where
\begin{equation*}
\tilde{\gamma}_k^2=s^2+2k\nu s\cos(\theta_1-\theta_2)+k^2\nu^2.
\end{equation*}
The dispersion relation \eqref{disprelintro} can be written as
\begin{equation}\label{disprel2intro}
\l_1^2\left(\frac{\rho}{\tanh(\sqrt{l_1^2+l_2^2})}+\frac{1}{\tanh(h\sqrt{l_1^2+l_2^2})}\right)-\left(\alpha+\beta(l_1^2+l_2^2)\right)\sqrt{l_1^2+l_2^2}=0,
\end{equation}
where
\begin{align}
l_1&=\nu k\cos(\theta_2)+s\cos(\theta_1),\label{l1intro}\\
l_2&=\nu k\sin(\theta_2)+s\sin(\theta_1)\label{l2intro}.
\end{align}
\begin{figure}\centering
\includegraphics{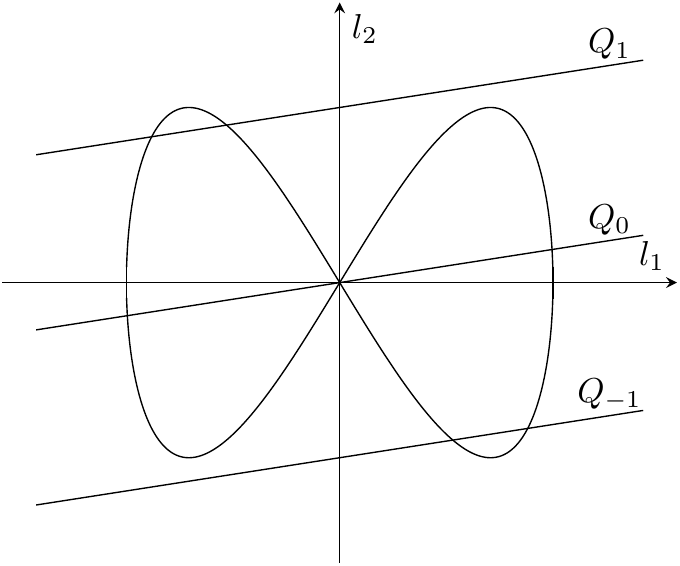}
\caption{Intersections between the real branches of $C_{dr}$ and the lines $Q_{\pm 1}$, $Q_0$.}
\label{inter2}
\end{figure}
The solution set of \eqref{disprelintro} can therefore be interpreted geometrically as in \cite{HaragusCourcelle2000}, namely that an imaginary number $\i s$ is a mode $k$ eigenvalue if and only if the line
\begin{equation*}
Q_k=\{(l_1,l_2)\ :\ l_1=\nu k\cos(\theta_2)+s\cos(\theta_1),\ l_2=\nu k\sin(\theta_2)+s\sin(\theta_1),\ s\in \mathbb{R}\},
\end{equation*}
intersects the real solution branch $C_{dr}$ of \eqref{disprel2intro}, see Figure \ref{inter2}. Due to this it is possible to obtain the same bifurcation scenarios in the internal wave setting as in the surface wave setting and we refer to \cite{Groves2003} where they list all the possible bifurcation scenarios for surface waves, involving mode $0$ and mode $\pm 1$ eigenvalues. However in the present work we focus on two particular cases.
One of these cases is a Hamiltonian-Hopf bifurcation involving mode $\pm 1$ eigenvalues. The bifurcation is achieved by choosing $\nu$ such that $\pm \i s$ are algebraically double mode $\pm 1$ eigenvalues. This occurs precisely when the lines $Q_1$,$Q_{-1}$ are tangential to $C_{dr}$ (see Figure \ref{inter1}). To see that this is a Hamiltonian-Hopf bifurcation, we first note that $\nu$ determines where $Q_k$ intersects the $l_2$-axis. In particular, when $\nu$ is large enough the lines $Q_k$ will not intersect $C_{dr}$ for $\abs{k}\geq 1$. So in particular there are no mode $\pm 1$ eigenvalues for such values of $\nu$, instead there is a plus minus complex conjugate quartet of complex mode $\pm 1$ eigenvalues. When $\nu$ is decreased there will be some critical value $\nu_0$ so that $Q_1, Q_{-1}$ are tangential to $C_{dr}$ which yields the algebraically double mode $\pm 1$ eigenvalues $\pm \i s$. When $\nu$ is decreased further the lines $Q_{1},Q_{-1}$ intersect $C_{dr}$ in two points each (see Figure \ref{inter2}), which means that there are four algebraically simple mode $\pm 1$ eigenvalues $\pm \i s_1,\pm \i s_2$. For equations that do not have a Hamiltonian structure, this bifurcation is called a $(\i s)^2$ resonance. The Hamilton-Hopf bifurcation is illustrated in \ref{changeinspectrumhamhopf}.
\begin{figure}
\centering
\subfigure[$\nu>\nu_0$]{
\includegraphics[scale=0.6]{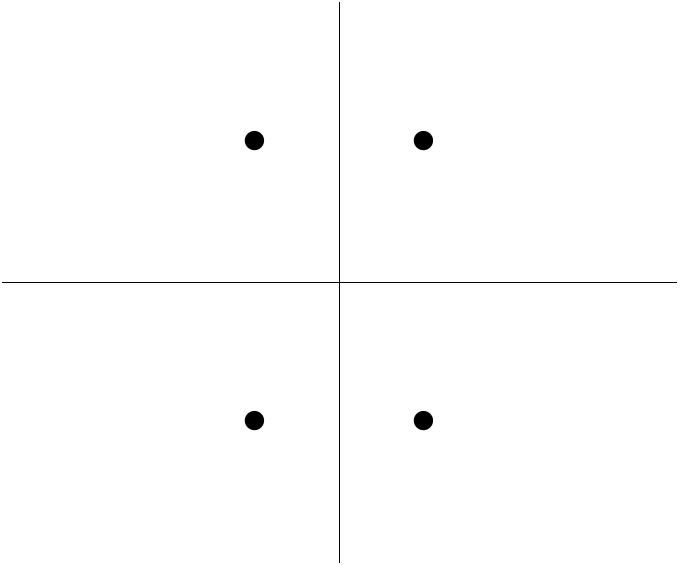}
}
\subfigure[$\nu=\nu_0$]{
\includegraphics[scale=0.6]{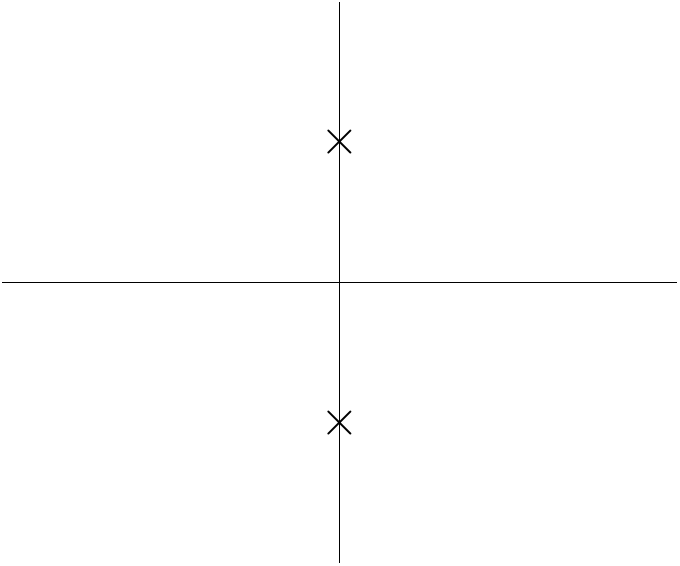}
}
\subfigure[$\nu<\nu_0$]{
\includegraphics[scale=0.6]{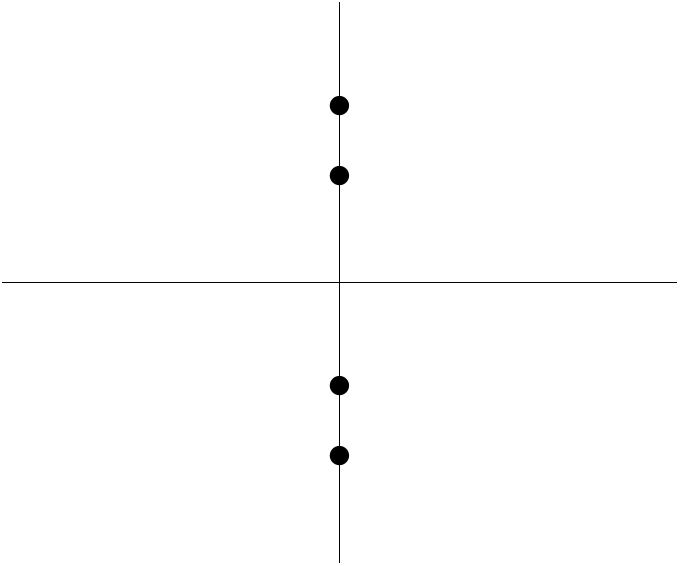}
}
\label{changeinspectrumhamhopf}
\caption{The Hamiltonian-Hopf bifurcation.}
\end{figure} 
The other bifurcation scenario we will consider is the following type of $00(\i s_2)$ resonance. Assume first that the lines $Q_1$ and $Q_{-1}$ intersect $C_{dr}$ in two distinct points each, so that $K$ has the mode $\pm 1$ eigenvalues $\pm \i s_1, \pm \i s_2$. Assume next that there is a critical value $\nu_0$ of $\nu$ such that $s_1=0$. Then $0$ is a mode $\pm 1$ eigenvalue so in particular it is of geometric multiplicity $2$. When $\nu$ is decreased through this critical value $\nu_0$, the following change in the spectrum of $K$ occurs. For $\nu>\nu_0$, $K$ has the mode $\pm 1$ eigenvalues $\pm \i s_1$. When $\nu$ is decreased to $\nu_0$ the eigenvalues $\pm \i s_1$ collide at the origin and form the geometrically double eigenvalue $0$. When $\nu$ is decreased further $K$ will again have the eigenvalues $\pm \i s_1$, however now $\i s_1$ is a mode $-1$ eigenvalue and $-\i s_1$ is a mode $1$ eigenvalue. The notation is changed in a natural way depending on the spectrum of $K$. For example, if $K$ has another pair of mode $\pm k$ eigenvalues $\pm \i s_3$ we denote the resonance by $00(\i s_2)(\i s_3)$. We illustrate the $00(\i s_2)$ resonance in Figure \ref{changeinspectrumdoublyper}.
\begin{figure}
\centering
\subfigure[$\nu>\nu_0$]{
\includegraphics[scale=0.6]{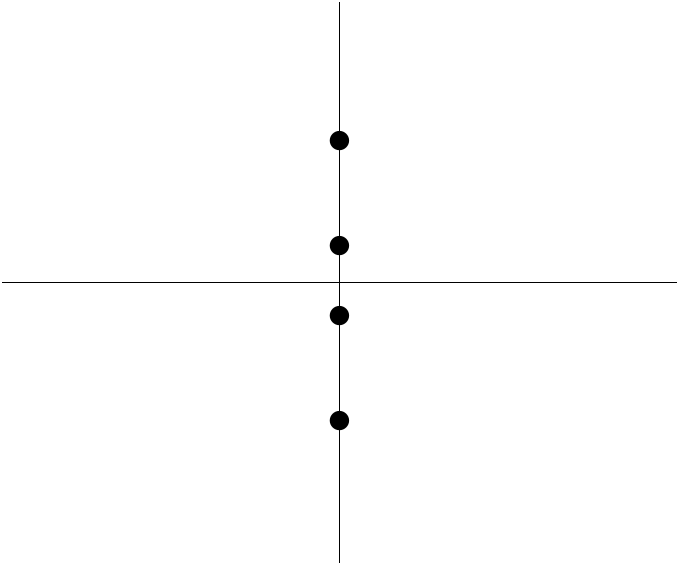}
}
\subfigure[$\nu=\nu_0$]{
\includegraphics[scale=0.6]{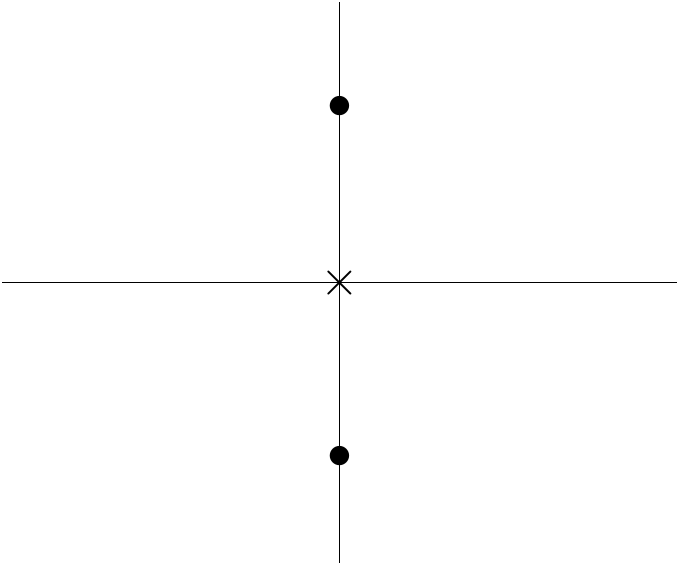}
}
\subfigure[$\nu<\nu_0$]{
\includegraphics[scale=0.6]{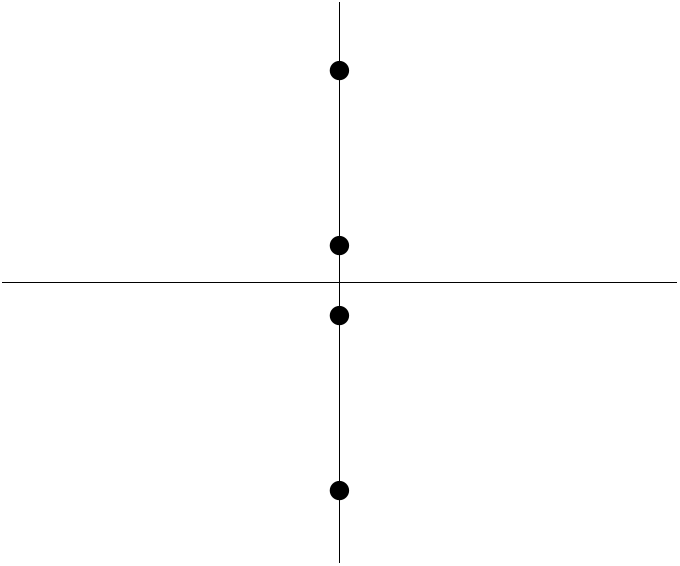}
}
\label{changeinspectrumdoublyper}
\caption{The $00(\mathrm{i}s_2)$ resonance.}
\end{figure} 

Section \ref{cmtsection} consists of a statement of the center-manifold theorem, and verification of the hypotheses of the theorem. We are using a version of the theorem which is due to Mielke \cite{MI}. In particular the theorem preserves the Hamiltonian structure of equation \eqref{modelhameq} so that the finite-dimensional reduced system also has a Hamiltonian structure.
\begin{figure}\centering
\includegraphics[scale=0.9]{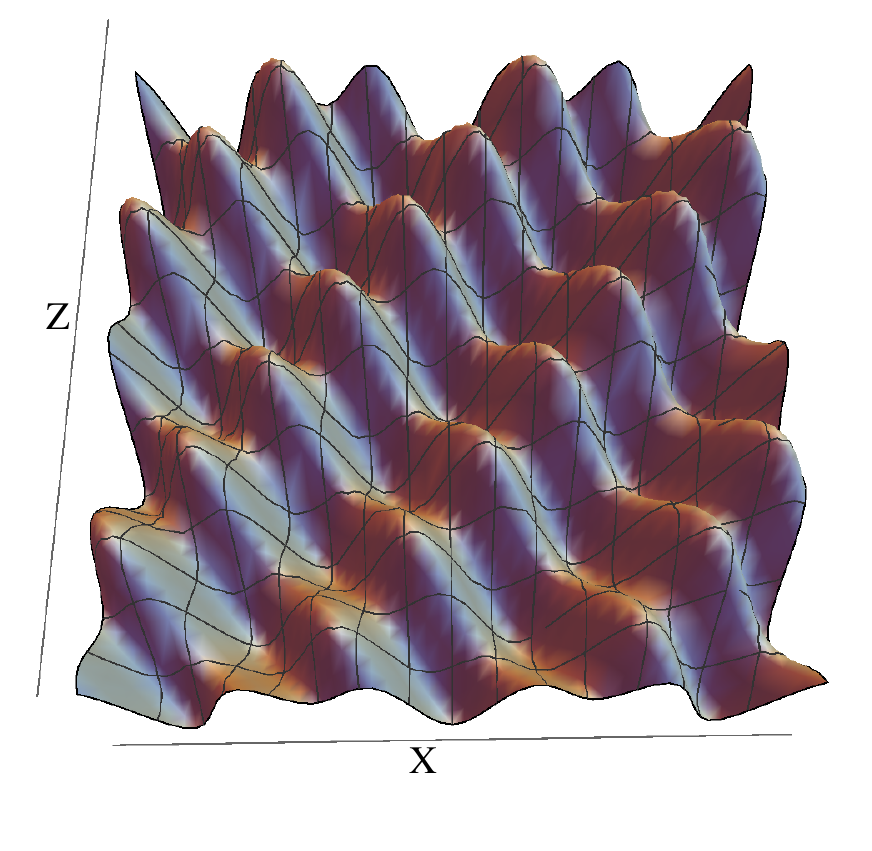}
\caption{Example of a doubly periodic wave.}
\label{doubleper}
\end{figure} 
In section \ref{periodic solutions} we construct doubly periodic internal waves, that is waves that are periodic in both $x$ and $z$, see Figure \ref{doubleper}. We do this by considering a $00(\i s)(\i \kappa_0)$ resonance, where $\pm \i s$ are mode $\pm 1$ eigenvalues and $\pm \i \kappa_0$ are mode $0$ eigenvalues, all algebraically simple. We also need to assume that $\i s$ is nonresonant with $\i \kappa_0$, that is $s\neq m\kappa_0$ for all $m\in \mathbb{Z}$. Note that the Lyapunov center theorem (see for example \cite{Ambrosetti}) cannot be applied here since this requires all eigenvalues to be nonresonant. This is clearly violated in our case due to the eigenvalue $0$. We approach this problem by first performing a center-manifold reduction, which gives us a finite-dimensional Hamiltonian system. This system is further reduced by a variational Lyapunov-Schmidt reduction. The existence of solutions is then established by an application of the implicit function theorem. 
The method employed in the present paper for proving existence of doubly periodic waves is different from the one used in \cite{Groves2003}. In that paper the authors apply the Lyapunov center theorem in order to find such waves. This is done by fixing the parameters so that $K$ has the mode $\pm 1$ eigenvalues $\pm \i \kappa$ which are nonresonant with all other eigenvalues.  After performing a center-manifold reduction, the Lyapunov center theorem yields solutions which are periodic in $x$ and $z$, with periods respectively near $2\pi/\kappa$ and equal to $2\pi/\nu$. At the linear level the solutions are given by linear combinations of
\begin{equation}\label{ghsolutions}
\exp(\pm \i (\kappa x+\nu z)).
\end{equation}
These yield solutions of the nonlinear problem that depend on both $x$ and $z$, but only through the combination $\kappa x+\nu z$. We note that
\begin{equation}\label{l1l2}
\kappa x+\nu z=l_1X+l_2Z,
\end{equation}
and there exist an angle $\tilde{\theta}_1$ such that 
\begin{equation*}
(l_1,l_2)=\sqrt{l_1^2+l_2^2}(\cos(\tilde{\theta}_1),\sin(\tilde{\theta}_1))=\tilde{\gamma}_1(\cos(\tilde{\theta}_1),\sin(\tilde{\theta}_1)).
\end{equation*}
From \eqref{l1l2} we then obtain
\begin{equation*}
\kappa x+\nu z=\tilde{\gamma}_1(\cos(\tilde{\theta}_1)X+\sin(\tilde{\theta}_1)Z)=:\tilde{\gamma}_1\tilde{x}.
\end{equation*}
A calculation shows that these waves are solutions of \eqref{goveq1}--\eqref{goveq8} that depend upon $\tilde{x}$, $Y$. Hence, these solutions are oblique line waves. In comparison, the solutions at the linear level obtained in the present paper are linear combinations of
\begin{equation*}
\exp(\pm \i \kappa_0x),\ \exp(\pm \i \nu_0z),
\end{equation*}
which yield solutions of the nonlinear problem that are doubly periodic in $x$ and $z$ and genuinely three-dimensional (see Theorem \ref{doublyperiodic_thm}).

In section \ref{Hamiltonian-Hopf bifurcation} we consider the Hamiltonian-Hopf bifurcation. After performing a center-manifold reduction and applying normal form theory, we obtain the reduced Hamiltonian system
\begin{align}
A_x&=\i sA+B+\i c_3^1\mu A+d_2^0\i A\abs{A}^2-2d_3^0A(A\bar{B}-\bar{A}B)+\mathcal{O}(\abs{(A,B)}\abs{(A,B,\mu)}^3)\label{introexphameqA},\\
B_x&=\i sB+\i c_3\mu B-c_2^1\mu A-2d_1^0A\abs{A}^2-\i d_2^0A^2\bar{B}+2d_2^0B\abs{A}^2-2d_3^0B(A\bar{B}-\bar{A}B)\nonumber\\
&\quad +\mathcal{O}(\abs{(A,B)}\abs{(A,B,\mu)}^3),\label{introexphameqB}
\end{align}
where $A,B$ are coordinates on the center manifold and $\mu$ is a bifurcation parameter. The solution set of \eqref{introexphameqA}--\eqref{introexphameqB} depends upon the signs of the coefficients $c_2^1$ and $d_1^0$. We find that both of the cases $c_2^1<0$, $d_1^0>0$ and $c_2^1<0$, $d_1^0<0$ can occur, depending on the the different parameters involved. So the situation here is analogous to the two-dimensional case studied in \cite{nilsson2016}. When $c_2^1<0$, $d_1^0>0$, $\mu>0$ the system \eqref{introexphameqA}--\eqref{introexphameqB} has two bright solitary wave solutions that each generate a one-parameter family of multipulse solutions, and when $c_2^1<0$, $d_1^0<0$, $\mu<0$ the system \eqref{introexphameqA}--\eqref{introexphameqB} has a one-parameter family of dark solitary waves. See Figure \ref{darkbrightwaveprofiles} for sketches of the different types of solutions.
\section{Spatial dynamics formulation of the travelling water wave problem}\label{spatialdyn}
Introduce in \eqref{goveq1}--\eqref{goveq8} the non-dimensional variables
\begin{equation*}
(X',Y',Z')=\frac{1}{h_1}(X,Y,Z),\ \eta'(X',Z')=\frac{1}{h_1}\eta(X,Z),\ \phi_i'(X',Y',Z')=\frac{1}{h_1c}\phi_i(X,Y,Z).
\end{equation*}
This gives us the system of equations,
\begin{align*}
\Delta\phi_1&=0,\quad \text{for }\eta<Y<1,\\
\Delta\phi_2&=0,\quad \text{for }-h<Y<\eta,
\end{align*}
with boundary conditions
\begin{align*}
&\phi_{1Y}=0&&\text{on }Y=1,\\
&\phi_{2Y}=0&&\text{on }Y=-h,\\
&\phi_{1Y}=-\eta_X+\eta_X\phi_{1X}+\eta_Z\phi_{1Z}&&\text{on }Y=\eta,\\
&\phi_{2Y}=-\eta_X+\eta_X\phi_{2X}+\eta_Z\phi_{2Z}&&\text{on }Y=\eta,\\
&-\phi_{2X}+\frac{1}{2}\abs{\nabla\phi_2}^2-\rho\left(-\phi_{1X}+\frac{1}{2}\abs{\nabla\phi_1}^2\right)+\alpha\eta=\beta\left(\frac{\eta_X}{\sqrt{1+\eta_X^2+\eta_Z^2}}\right)_X\\
&\quad+\beta\left(\frac{\eta_Z}{\sqrt{1+\eta_X^2+\eta_Z^2}}\right)_Z&&\text{on }Y=\eta,
\end{align*}
where $h=h_2/h_1$, $\rho=\rho_1/\rho_2$, $\alpha=gh_1(1-\rho)/c^2$, $\beta=\sigma/(h_1\rho_2c^2)$ and we have dropped the prime for notational simplicity. Next we introduce the coordinate system $(x,z)$ given in \eqref{newcoordinates} and look for solutions of the form
\begin{equation*}
\tilde{\eta}(x,z)=\eta(X,Z),\  \tilde{\phi}_i(x,y,z)=\phi_i(X,Y,Z),
\end{equation*}
such that $\tilde{\eta}$ and $\tilde{\phi}_i$ are periodic in $z$, with period $P_z$. We also obtain a fixed domain by defining
\begin{equation*}
y(x,z)=\begin{cases}
\frac{Y-1}{\eta-1}\ &\eta<Y<1,\\
\frac{Y+h}{\eta+h}\ &-h<Y<\eta.
\end{cases}
\end{equation*}
Finally, the period is normalized to be $2\pi$ and if we let $\nu=2\pi/P_z$ the governing equations become (with the tilde removed)
\begin{align}
&\phi_{1xx}-\frac{2\eta_xy\phi_{1yx}}{\eta-1}+\left(-\frac{\eta_{xx}y}{\eta-1}+\frac{2\eta_x^2y}{(\eta-1)^2}\right)\phi_{1y}+\frac{\eta_x^2y^2\phi_{1yy}}{(\eta-1)^2}+\frac{\phi_{1yy}}{(\eta-1)^2}+\nonumber \\
&\quad+\nu^2\bigg[\phi_{1zz}-\frac{2\eta_zy\phi_{1yz}}{\eta-1}+\left(-\frac{\eta_{zz}y}{\eta-1}+\frac{2\eta_z^2y}{(\eta-1)^2}\right)\phi_{1y}+\frac{\eta_z^2y^2\phi_{1yy}}{(\eta-1)^2}\bigg]\nonumber\\
&\quad+2\nu\cos(\theta_1-\theta_2)\bigg[\phi_{1xz}-\frac{\eta_xy\phi_{1yz}}{\eta-1}+\left(-\frac{\eta_{zx}y}{\eta-1}+\frac{2\eta_{x}\eta_zy}{(\eta-1)^2}\right)\phi_{1y}-\frac{\eta_zy\phi_{1xy}}{\eta-1}\nonumber\\
&\quad+\frac{\eta_z\eta_xy^2\phi_{1yy}}{(\eta-1)^2}\bigg]=0\quad \text{for }0<y<1,\label{goveqfinal1}\\
&\phi_{2xx}-\frac{2\eta_xy\phi_{2yx}}{\eta+h}+\left(-\frac{\eta_{xx}y}{\eta+h}+\frac{2\eta_x^2y}{(\eta+h)^2}\right)\phi_{2y}+\frac{\eta_x^2y^2\phi_{2yy}}{(\eta+h)^2}+\frac{\phi_{2yy}}{(\eta+h)^2}\nonumber\\ 
&\quad+\nu^2\bigg[\phi_{2zz}-\frac{2\eta_zy\phi_{2yz}}{\eta+h}+\left(-\frac{\eta_{zz}y}{\eta+h}+\frac{2\eta_z^2y}{(\eta+h)^2}\right)\phi_{2y}+\frac{\eta_z^2y^2\phi_{2yy}}{(\eta+h)^2}\bigg]\nonumber\\
&\quad+2\nu\cos(\theta_1-\theta_2)\bigg[\phi_{2xz}-\frac{\eta_xy\phi_{2yz}}{\eta+h}+\left(-\frac{\eta_{zx}y}{\eta+h}+\frac{2\eta_{x}\eta_zy}{(\eta+h)^2}\right)\phi_{2y}-\frac{\eta_zy\phi_{2xy}}{\eta+h}\nonumber\\
&\quad+\frac{\eta_z\eta_xy^2\phi_{2yy}}{(\eta+h)^2}\bigg]=0\quad\text{for }0<y<1,
\end{align}
with boundary conditions
\begin{align}
&\phi_{1y}=0&&\text{on }y=0,\\
&\phi_{2y}=0&&\text{on }y=0,\\
&\frac{\phi_{1y}}{\eta-1}=\eta_x\left(\phi_{1x}-\frac{\eta_xy\phi_{1y}}{\eta-1}\right)+\nu^2\eta_z\left(\phi_{1z}-\frac{\eta_zy\phi_{1y}}{\eta-1}\right)-\cos(\theta_1)\eta_x\nonumber\\
& -\nu\cos(\theta_2)\eta_z+ \nu\cos(\theta_1-\theta_2)\left(\eta_x\left(\phi_{1z}-\frac{\eta_zy\phi_{1y}}{\eta-1}\right)+\eta_z\left(\phi_{1x}-\frac{\eta_xy\phi_{1y}}{\eta-1}\right)\right)&&\text{on }y=1,\\
&\frac{\phi_{2y}}{\eta+h}=\eta_x\left(\phi_{2x}-\frac{\eta_xy\phi_{2y}}{\eta+h}\right)+\nu^2\eta_z\left(\phi_{2z}-\frac{\eta_zy\phi_{2y}}{\eta+h}\right)-\cos(\theta_1)\eta_x\nonumber\\
& -\nu\cos(\theta_2)\eta_z+ \nu\cos(\theta_1-\theta_2)\left(\eta_x\left(\phi_{2z}-\frac{\eta_zy\phi_{2y}}{\eta+h}\right)+\eta_z\left(\phi_{2x}-\frac{\eta_xy\phi_{2y}}{\eta+h}\right)\right)&&\text{on }y=1,
\end{align}
\begin{align}
&-\cos(\theta_1)\left(\phi_{2x}-\frac{\eta_xy\phi_{2y}}{\eta+h}\right)-\nu\cos(\theta_2)\left(\phi_{2z}-\frac{\eta_zy\phi_{2y}}{\eta+h}\right)+\frac{1}{2}\bigg[\left(\phi_{2x}-\frac{\eta_xy\phi_{2y}}{\eta+h}\right)^2\nonumber\\
&\quad+\frac{\phi_{2y}^2}{(\eta+h)^2}+\nu^2\left(\phi_{2z}-\frac{\eta_zy\phi_{2y}}{\eta+h}\right)^2+2\nu\cos(\theta_1-\theta_2)\left(\phi_{2x}-\frac{\eta_xy\phi_{2y}}{\eta+h}\right)\left(\phi_{2z}-\frac{\eta_zy\phi_{2y}}{\eta+h}\right)\bigg]\nonumber\\
&\quad-\rho\Bigg\{-\cos(\theta_1)\left(\phi_{1x}-\frac{\eta_xy\phi_{1y}}{\eta-1}\right)-\nu\cos(\theta_2)\left(\phi_{1z}-\frac{\eta_zy\phi_{1y}}{\eta-1}\right)+\frac{1}{2}\bigg[\left(\phi_{1x}-\frac{\eta_xy\phi_{1y}}{\eta-1}\right)^2\nonumber\\
&\quad+\frac{\phi_{1y}^2}{(\eta-1)^2}+\nu^2\left(\phi_{1z}-\frac{\eta_zy\phi_{1y}}{\eta-1}\right)^2+2\nu\cos(\theta_1-\theta_2)\left(\phi_{1x}-\frac{\eta_xy\phi_{1y}}{\eta-1}\right)\left(\phi_{1z}-\frac{\eta_zy\phi_{1y}}{\eta-1}\right)\bigg]\Bigg\}\nonumber\\
&\quad +\alpha\eta=\beta\left[\left(\frac{\eta_x}{Q}\right)_x+\nu^2\left(\frac{\eta_z}{Q}\right)_z+\nu\cos(\theta_1-\theta_2)\left(\left(\frac{\eta_z}{Q}\right)_x+\left(\frac{\eta_x}{Q}\right)_z\right)\right]\quad\quad\ \text{on }y=1,\label{goveqfinal8}
\end{align}
where  $Q=\sqrt{1+\eta_x^2+\nu^2\eta_z^2+2\nu\cos(\theta_1-\theta_2)\eta_x\eta_z}$.

The energy and momentum associated with this system are given by
\begin{align*}
E&=\frac{\rho_1h_1^3c^2}{2\nu}\int_\mathbb{R}\int_0^{2\pi}\int_0^1\bigg[\left(\phi_{1x}-\frac{\eta_xy\phi_{1y}}{\eta-1}\right)^2+\frac{\phi_{1y}^2}{(\eta-1)^2}+\nu^2\left(\phi_{1z}-\frac{\eta_y\phi_{1y}}{\eta-1}\right)^2\\
&\quad + 2\nu\cos(\theta_1-\theta_2)\left(\phi_{1x}-\frac{\eta_xy\phi_{1y}}{\eta-1}\right)\left(\phi_{1z}-\frac{\eta_zy\phi_{1y}}{\eta-1}\right)\bigg](1-\eta)\ \mathrm{d}y\ \mathrm{d}z\ \mathrm{d}x\\
&\quad +\frac{\rho_2c^2h_1^3}{2\nu}\int_\mathbb{R}\int_0^{2\pi}\int_0^1\bigg[\left(\phi_{2x}-\frac{\eta_xy\phi_{2y}}{\eta+h}\right)^2+\frac{\phi_{2y}^2}{(\eta+h)^2}+\nu^2\left(\phi_{2z}-\frac{\eta_y\phi_{2y}}{\eta+h}\right)^2\\
&\quad +2\nu\cos(\theta_1-\theta_2)\left(\phi_{2x}-\frac{\eta_xy\phi_{2y}}{\eta+h}\right)\left(\phi_{2z}-\frac{\eta_zy\phi_{2y}}{\eta+h}\right)\bigg](\eta+h)\ \mathrm{d}y\ \mathrm{d}z\ \mathrm{d}x\\
&\quad +\frac{g(\rho_2-\rho_1)h_1^4}{2\nu}\int_\mathbb{R}\int_0^{2\pi}\eta^2\ \mathrm{d}z\ \mathrm{d}x+\frac{\sigma h_1^2}{\nu}\int_\mathbb{R}\int_0^{2\pi}Q-1\ \mathrm{d}z \ \mathrm{d}x,
\end{align*}
\begin{align*}
P&=\frac{\rho_1h_1^3c}{\nu}\int_\mathbb{R}\int_0^{2\pi}\int_0^1\bigg[\cos(\theta_1)\left(\phi_{1x}-\frac{\eta_xy\phi_{1y}}{\eta-1}\right)+\nu\cos(\theta_1)\left(\phi_{1z}-\frac{\eta_zy\phi_{1y}}{\eta-1}\right)\bigg](1-\eta)\ \mathrm{d}y\ \mathrm{d}z \ \mathrm{d}x\\
&\quad+\frac{\rho_2h_1^3c}{\nu}\int_\mathbb{R}\int_0^{2\pi}\int_0^1\bigg[\cos(\theta_1)\left(\phi_{2x}-\frac{\eta_xy\phi_{2y}}{\eta+h}\right)+\nu\cos(\theta_2)\left(\phi_{1z}-\frac{\eta_zy\phi_{2y}}{\eta+h}\right)\bigg](\eta+h)\ \mathrm{d}y\ \mathrm{d}z \ \mathrm{d}x.
\end{align*}
The solutions we are interested in are critical points of the functional $E-cP$. This is an action integral, with Lagrangian
\begin{align*}
L(\eta,\eta_x,\phi_1,\phi_{1x},\phi_2,\phi_{2x})&=\int_0^{2\pi}\int_0^1\frac{\rho}{2}\bigg[\left(\phi_{1x}-\frac{\eta_xy\phi_{1y}}{\eta-1}-\cos(\theta_1)\right)^2+\left(\nu\left(\phi_{1z}-\frac{\eta_zy\phi_{1y}}{\eta-1}\right)-\cos(\theta_2)\right)^2\\
&\quad+\frac{\phi_{1y}^2}{(\eta-1)^2}+2\nu\cos(\theta_1-\theta_2)\left(\phi_{1x}-\frac{\eta_xy\phi_{1y}}{\eta-1}\right)\left(\phi_{1z}-\frac{\eta_zy\phi_{1y}}{\eta-1}\right)\\
&\quad-\cos^2(\theta_1)-\cos^2(\theta_2)\bigg](1-\eta)\ \mathrm{d}y \ \mathrm{d}z +\int_0^{2\pi}\int_0^1\frac{1}{2}\bigg[\left(\phi_{2x}-\frac{\eta_xy\phi_{2y}}{\eta+h}-\cos(\theta_1)\right)^2\\
&\quad +\left(\nu\left(\phi_{2z}-\frac{\eta_zy\phi_{2y}}{\eta+h}\right)-\cos(\theta_2)\right)^2+\frac{\phi_{2y}^2}{(\eta+h)^2}\\
&\quad+2\nu\cos(\theta_1-\theta_2)\left(\phi_{2x}-\frac{\eta_xy\phi_{2y}}{\eta+h}\right)\left(\phi_{2z}-\frac{\eta_zy\phi_{2y}}{\eta+h}\right)\\
&\quad-\cos^2(\theta_1)-\cos^2(\theta_2)\bigg](\eta+h)\ \mathrm{d}y \ \mathrm{d}z+\frac{\alpha}{2}\int_0^{2\pi}\eta^2\ \mathrm{d}z +\beta\int_0^{2\pi}Q-1\ \mathrm{d}z.
\end{align*}
A Hamiltonian formulation of \eqref{goveqfinal1}--\eqref{goveqfinal8} is obtained via the Legendre transform
\begin{align*}
\psi_1&:=\frac{\delta L}{\delta \phi_{1x}}=\rho\bigg[\phi_{1x}-\frac{\eta_xy\phi_{1y}}{\eta-1}-\cos(\theta_1)+\nu\cos(\theta_1-\theta_2)\left(\phi_{1z}-\frac{\eta_zy\phi_{1y}}{\eta-1}\right)\bigg](1-\eta),\\
\psi_2&:=\frac{\delta L}{\delta \phi_{2x}}=\bigg[\phi_{2x}-\frac{\eta_xy\phi_{2y}}{\eta+h}-\cos(\theta_1)+\nu\cos(\theta_1-\theta_2)\left(\phi_{2z}-\frac{\eta_zy\phi_{2y}}{\eta+h}\right)\bigg](\eta+h),\\
\omega&:=\frac{\delta L}{\delta \eta_x}=-\int_0^1\frac{y\phi_{1y}\psi_1}{\eta-1}\ \mathrm{d}y -\int_0^1\frac{y\phi_{2y}\psi_2}{\eta+h}\ \mathrm{d}y+\frac{\beta(\eta_x+\nu\cos(\theta_1-\theta_2)\eta_z)}{Q}.
\end{align*}
The Hamiltonian $\mathcal{H}$ is then defined by
\begin{align*}
\mathcal{H}(\eta,\omega,\phi_1,\psi_1,\phi_2,\psi_2)&=\int_0^{2\pi}\int_0^1\psi_1\phi_{1x}\ \mathrm{d}y\ \mathrm{d}z+\int_0^{2\pi}\int_0^1\psi_2\phi_{2x}\ \mathrm{d}y\ \mathrm{d}z+\int_0^{2\pi}\omega\eta_x\ \mathrm{d}z\\
&\quad -L(\eta,\eta_x,\phi_1,\phi_{1x},\phi_2,\phi_{2x})\\
&=\int_0^{2\pi}\int_0^1 \bigg[\frac{1}{2\rho(1-\eta)}\left(\psi_1-\rho(1-\eta)\nu\cos(\theta_1-\theta_2)\left(\phi_{1z}-\frac{\eta_zy\phi_{1y}}{\eta-1}\right)\right)^2\\
&\quad+\psi_1\cos(\theta_1)-\frac{\rho(1-\eta)}{2}\left(\nu\left(\phi_{1z}-\frac{\eta_zy\phi_{1y}}{\eta-1}\right)-\cos(\theta_2)\right)^2-\frac{\rho\phi_{1y}^2}{2(1-\eta)}\\
&\quad-\rho\nu\cos(\theta_1-\theta_2)(1-\eta)\cos(\theta_1)\left(\phi_{1z}-\frac{\eta_zy\phi_{1y}}{\eta-1}\right)\\
&\quad+\frac{\rho(1-\eta)}{2}\left(\cos^2(\theta_1)+\cos^2(\theta_2)\right)\bigg]\ \mathrm{d}y \ \mathrm{d}z\\
&\quad +\int_0^{2\pi}\int_0^1 \bigg[\frac{1}{2(\eta+h)}\left(\psi_2-(\eta+h)\nu\cos(\theta_1-\theta_2)\left(\phi_{1z}-\frac{\eta_zy\phi_{1y}}{\eta+h}\right)\right)^2\\
&\quad+\psi_2\cos(\theta_1)-\frac{(\eta+h)}{2}\left(\nu\left(\phi_{2z}-\frac{\eta_zy\phi_{2y}}{\eta+h}\right)-\cos(\theta_2)\right)^2-\frac{\phi_{2y}^2}{2(\eta+h)}\\
&\quad-\nu\cos(\theta_1-\theta_2)(\eta+h)\cos(\theta_1)\left(\phi_{2z}-\frac{\eta_zy\phi_{2y}}{\eta+h}\right)\\
&\quad+\frac{(\eta+h)}{2}\left(\cos^2(\theta_1)+\cos^2(\theta_2)\right)\bigg]\ \mathrm{d}y \ \mathrm{d}z -\frac{\alpha}{2}\int_0^{2\pi}\eta^2\ dz-\int_0^{2\pi}\bigg[\nu\cos(\theta_1-\theta_2)\eta_z\bar{\omega}\\
&\quad+\sqrt{\beta^2-\bar{\omega}^2}\sqrt{\nu^2\sin^2(\theta_1-\theta_2)\eta_z^2+1}-\beta\ \bigg]\ \mathrm{d}z,
\end{align*}
where
\begin{equation*}
\bar{\omega}=\omega+\int_0^1\frac{y\phi_{1y}\psi_1}{\eta-1}\ \mathrm{d}y +\int_0^1\frac{y\phi_{2y}\psi_2}{\eta+h}\ \mathrm{d}y.
\end{equation*}
For $s\geq 0$, define
\begin{equation*}
X_s=H_{\text{per}}^{s+1}(S)\times H_{\text{per}}^s(S)\times H_{\text{per}}^{s+1}(\Sigma)\times H_{\text{per}}^{s}(\Sigma)\times H_{\text{per}}^{s+1}(\Sigma)\times H_{\text{per}}^{s}(\Sigma),
\end{equation*}
where $S=(0,2\pi)$, $\Sigma=(0,2\pi)\times (0,1)$ and
\begin{align*}
H_{\text{per}}^s(S)&=\{f\in H_\text{loc}^s(\mathbb{R})\ :\ f(z+2\pi)=f(z), \ z\in\mathbb{R}\},\\
H_\text{per}^s(\Sigma)&=\{f\in H_\text{loc}^s((0,1)\times\mathbb{R})\ :\ f(y,z+2\pi)=f(y,z),\ y\in (0,1),\ z\in\mathbb{R}\}.
\end{align*}
Let $M=X_0$, $m\in M$ and let $v=(\eta,\omega,\phi_1,\psi_1,\phi_2,\psi_2)\in T_mM$. On $T_mM\times T_mM$ we define the position-independent symplectic form
\begin{align}
\Omega(v,v^*)&=\int_S(\omega^*\eta-\eta^*\omega)\ \mathrm{d}z +\int_\Sigma(\psi_1^*\phi_1-\phi_1^*\psi_1)\ \mathrm{d}y \ \mathrm{d}z+\int_\Sigma(\psi_2^*\phi_2-\phi_2^*\psi_2)\ \mathrm{d}y\ \mathrm{d}z\label{omega1}.
\end{align}
As in \cite{Bagri2014} we observe that $(M,\Omega)$ is a symplectic manifold and that the set
\begin{equation*}
N=\{m\in X_1\text{ }:\text{ } \vert \bar{\omega}(z)\vert<\beta,\text{ } -h<\eta(z)<1\}
\end{equation*}
is a manifold domain of $M$ with $H\in C^\infty(N,\mathbb{R})$. The triple $(M,H,\Omega)$ is therefore a Hamiltonian system. Note that in for example the papers \cite{Groves2003,Groves2007} they use $X_s,X_{s+1}$, for some $s\in(0,1/2)$, to construct the symplectic manifold. However, it was shown in \cite{Bagri2014} that it is possible to use the spaces $X_0,X_1$ and still obtain a well defined Hamiltonian system. The Hamiltonian vector field $v_\mathcal{H}$ and its domain $\mathcal{D}(v_\mathcal{H})$ is defined by
\begin{equation*}
\mathcal{D}(v_\mathcal{H}):=\{m\in N\ \vert\text{ } \exists (v_\mathcal{H})_m\in T_mM \text{ such that } \mathrm{d}\mathcal{H}[m](v_m^*)=\Omega((v_\mathcal{H})_m,v_m^*) \text{ } \forall v_m^*\in T_mM\},
\end{equation*}
and Hamilton's equation is given by
\begin{equation}\label{heq}
\dot{\gamma}(x)=(v_\mathcal{H})_{\gamma(x)}. 
\end{equation}
Before writing down \eqref{heq} explicitly we introduce the new variables $\tilde{\psi}_1=\psi_1+\rho\cos(\theta_1)$, $\tilde{\psi}_2=\psi_2+h\cos(\theta_1)$, so that $(0,0,0,0,0,0)$ is an equilibrium solution of the resulting Hamiltonian system. Suppressing the tilde, the Hamiltonian is then given by
\begin{align*}
&\mathcal{H}(\eta,\omega,\phi_1,\psi_1,\phi_2,\psi_2)\\
& =\int_0^{2\pi}\int_0^1 \bigg[\frac{1}{2\rho(1-\eta)}\left(\psi_1-\rho\cos(\theta_1)-\rho\nu\cos(\theta_1-\theta_2)(1-\eta)\left(\phi_{1z}-\frac{\eta_zy\phi_{1y}}{\eta-1}\right)\right)^2\\
&\quad+\psi_1\cos(\theta_1)-\rho\cos^2(\theta_1) -\frac{\rho(1-\eta)}{2}\left(\nu\left(\phi_{1z}-\frac{\eta_zy\phi_{1y}}{\eta-1}\right)-\cos(\theta_2)\right)^2-\frac{\rho\phi_{1y}^2}{2(1-\eta)}\\
&\quad-\rho\nu\cos(\theta_1-\theta_2)\cos(\theta_1)(1-\eta)\left(\phi_{1z}-\frac{\eta_zy\phi_{1y}}{\eta-1}\right)+\frac{\rho(1-\eta)}{2}\left(\cos^2(\theta_1)+\cos^2(\theta_2)\right)\bigg]\ \mathrm{d}y \ \mathrm{d}z\\
&\quad +\int_0^{2\pi}\int_0^1 \bigg[\frac{1}{2(\eta+h)}\left(\psi_2-h\cos(\theta_1)-\nu\cos(\theta_1-\theta_2)(\eta+h)\left(\phi_{1z}-\frac{\eta_zy\phi_{1y}}{\eta+h}\right)\right)^2\\
&\quad+\psi_2\cos(\theta_1)-h\cos^2(\theta_1) -\frac{(\eta+h)}{2}\left(\nu\left(\phi_{2z}-\frac{\eta_zy\phi_{2y}}{\eta+h}\right)-\cos(\theta_2)\right)^2-\frac{\phi_{2y}^2}{2(\eta+h)}\\
&\quad-\nu\cos(\theta_1-\theta_2)\cos(\theta_1)(\eta+h)\left(\phi_{2z}-\frac{\eta_zy\phi_{2y}}{\eta+h}\right)+\frac{(\eta+h)}{2}\left(\cos^2(\theta_1)+\cos^2(\theta_2)\right)\bigg]\ \mathrm{d}y \ \mathrm{d}z\\
&\quad -\frac{\alpha}{2}\int_0^{2\pi}\eta^2\ \mathrm{d}z-\int_0^{2\pi}\nu\cos(\theta_1-\theta_2)\eta_z\bar{\omega}+\sqrt{\beta^2-\bar{\omega}^2}\sqrt{\nu^2\sin^2(\theta_1-\theta_2)\eta_z^2+1}-\beta\ \mathrm{d}z,
\end{align*}
where
\begin{equation*}
\bar{\omega}=\omega+\int_0^1\frac{y\phi_{1y}(\psi_1-\rho\cos(\theta_1))}{\eta-1}\ \mathrm{d}y +\int_0^1\frac{y\phi_{2y}(\psi_2-h\cos(\theta_1))}{\eta+h}\ \mathrm{d}y,
\end{equation*}
and Hamilton's equations become
\begin{align}
\dot{\eta}&=\bar{\omega}\sqrt{\frac{\nu^2\sin^2(\theta_1-\theta_2)\eta_z^2+1}{\beta^2-\bar{\omega}^2}}-\nu\cos(\theta_1-\theta_2)\eta_z,\label{hamvec12}\\
\dot{\omega}&=\int_0^1\Bigg[-\frac{(\psi_1-\rho\cos(\theta_1))^2}{2\rho(1-\eta)^2}+\frac{\rho\phi_{1y}^2}{2(1-\eta)^2}-\frac{\rho\nu^2\sin^2(\theta_1-\theta_2)}{2}\left(\phi_{1z}^2-\frac{\eta_z^2y^2\phi_{1y}^2}{(\eta-1)^2}\right)\nonumber\\
&\quad +\bar{\omega}\sqrt{\frac{\nu^2\sin^2(\theta_1-\theta_2)\eta_z^2+1}{\beta^2-\bar{\omega}^2}}\left(\frac{y\phi_{1y}(\psi_1-\rho\cos(\theta_1))}{(\eta-1)^2}\right)\nonumber\\
&\quad-\rho\nu^2\sin^2(\theta_1-\theta_2)\left[y\phi_{1y}\left(\phi_{1z}-\frac{\eta_zy\phi_{1y}}{\eta-1}\right)\right]_z\Bigg]\ \mathrm{d}y\nonumber\\
&\quad -\rho\nu\left[(\cos(\theta_1-\theta_2)\cos(\theta_1)-\cos(\theta_2))\phi_{1z}\right]_{y=1}+\frac{\rho\cos^2(\theta_1)}{2}\nonumber\\
&\quad+\int_0^1\Bigg[\frac{(\psi_2-h\cos(\theta_1))^2}{2(\eta+h)^2}-\frac{\phi_{2y}^2}{2(\eta+h)^2}+\frac{\nu^2\sin^2(\theta_1-\theta_2)}{2}\left(\phi_{2z}^2-\frac{\eta_z^2y^2\phi_{2y}^2}{(\eta+h)^2}\right)\nonumber\\
&\quad +\bar{\omega}\sqrt{\frac{\nu^2\sin^2(\theta_1-\theta_2)\eta_z^2+1}{\beta^2-\bar{\omega}^2}}\left(\frac{y\phi_{2y}(\psi_2-h\cos(\theta_1))}{(\eta+h)^2}\right)\nonumber\\
&\quad+\nu^2\sin^2(\theta_1-\theta_2)\left[y\phi_{2y}\left(\phi_{2z}-\frac{\eta_zy\phi_{2y}}{\eta+h}\right)\right]_z\bigg]\ \mathrm{d}y\nonumber\\
&\quad +\nu\left[(\cos(\theta_1-\theta_2)\cos(\theta_1)-\cos(\theta_2))\phi_{2z}\right]_{y=1}-\frac{\cos^2(\theta_1)}{2}+\alpha\eta\nonumber\\
&\quad-\nu^2\sin^2(\theta_1-\theta_2)\left[\sqrt{\frac{\beta^2-\bar{\omega}^2}{\nu^2\sin^2(\theta_1-\theta_2)\eta_z^2+1}}\eta_z\right]_z-\nu\cos(\theta_1-\theta_2)\omega_z,\\
\dot{\phi}_1&=\frac{\psi_1-\rho\cos(\theta_1)}{\rho(1-\eta)}-\nu\cos(\theta_1-\theta_2)\phi_{1z}+\bar{\omega}\sqrt{\frac{\nu^2\sin^2(\theta_1-\theta_2)\eta_z^2+1}{\beta^2-\bar{\omega}^2}}\left(\frac{y\phi_{1y}}{\eta-1}\right)+\cos(\theta_1),\\
\dot{\psi}_1&=-\frac{\rho\phi_{1yy}}{1-\eta}-\nu\cos(\theta_1-\theta_2)\psi_{1z}+\rho\nu^2\sin^2(\theta_1-\theta_2)\left[(\eta-1)\left(\phi_{1z}-\frac{\eta_zy\phi_{1y}}{\eta-1}\right)\right]_z\nonumber\\
&\quad +\bar{\omega}\sqrt{\frac{\nu^2\sin^2(\theta_1-\theta_2)\eta_z^2+1}{\beta^2-\bar{\omega}^2}}\left[\frac{y(\psi_{1}-\rho\cos(\theta_1))}{\eta-1}\right]_y\nonumber\\
&\quad-\rho\nu^2\sin^2(\theta_1-\theta_2)\eta_z\left[y\left(\phi_{1z}-\frac{\eta_zy\phi_{1y}}{\eta-1}\right)\right]_y,\\
\dot{\phi}_2&=\frac{\psi_2-h\cos(\theta_1)}{\eta+h}-\nu\cos(\theta_1-\theta_2)\phi_{2z}+\bar{\omega}\sqrt{\frac{\nu^2\sin^2(\theta_1-\theta_2)\eta_z^2+1}{\beta^2-\bar{\omega}^2}}\left(\frac{y\phi_{2y}}{\eta+h}\right)+\cos(\theta_1),\\
\dot{\psi}_2&=-\frac{\phi_{2yy}}{\eta+h}-\nu\cos(\theta_1-\theta_2)\psi_{2z}-\nu^2\sin^2(\theta_1-\theta_2)\left[(\eta+h)\left(\phi_{2z}-\frac{\eta_zy\phi_{2y}}{\eta+h}\right)\right]_z\nonumber\\
&\quad+\bar{\omega}\sqrt{\frac{\nu^2\sin^2(\theta_1-\theta_2)\eta_z^2+1}{\beta^2-\bar{\omega}^2}}\left[\frac{y(\psi_{2}-h\cos(\theta_1))}{\eta+h}\right]_y\nonumber\\
&\quad+\nu^2\sin^2(\theta_1-\theta_2)\eta_z\left[y\left(\phi_{2z}-\frac{\eta_zy\phi_{2y}}{\eta+h}\right)\right]_y.\label{hamvec62}
\end{align}
The domain $\mathcal{D}(v_\mathcal{H})$ consists of elements in $(\eta,\omega,\phi_1,\psi_1,\phi_2,\psi_2)\in N$ such that
\begin{align}
&\begin{cases}
&\phi_{1y}=0, \qquad \qquad \qquad \qquad\qquad \qquad \qquad  \qquad \qquad \qquad \qquad \qquad \qquad \qquad \ y=0,\\
&-\frac{\rho\phi_{1y}}{1-\eta}-\rho\nu^2\sin^2(\theta_1-\theta_2)\eta_z\left(\phi_{1z}-\frac{\eta_z\phi_{1y}}{\eta-1}\right)-\rho\nu(\cos(\theta_1-\theta_2)\cos(\theta_1)-\cos(\theta_2))\eta_z\\
&+\bar{\omega}\sqrt{\frac{\nu^2\sin^2(\theta_1-\theta_2)\eta_z^2+1}{\beta^2-\bar{\omega}^2}}\left(\frac{\psi_1-\rho\cos(\theta_1)}{\eta-1}\right)=0, \qquad \qquad \qquad \qquad \qquad \qquad \qquad y=1,
\end{cases}\label{bdry11}\\
&\begin{cases}
&\phi_{2y}=0,\qquad \qquad \qquad \qquad\qquad \qquad \qquad  \qquad \qquad \qquad \qquad \qquad \qquad \qquad \ y=0,\\
&-\frac{\phi_{2y}}{\eta+h}+\nu^2\sin^2(\theta_1-\theta_2)\eta_z\left(\phi_{2z}-\frac{\eta_z\phi_{2y}}{\eta+h}\right)+\nu(\cos(\theta_1-\theta_2)\cos(\theta_1)-\cos(\theta_2))\eta_z\\
&+\bar{\omega}\sqrt{\frac{\nu^2\sin^2(\theta_1-\theta_2)\eta_z^2+1}{\beta^2-\bar{\omega}^2}}\left(\frac{\psi_2-h\cos(\theta_1)}{\eta+h}\right)=0,\qquad \qquad \qquad \qquad \qquad \qquad \qquad y=1.
\end{cases}\label{bdry21}
\end{align}
In the present paper we will only consider bifurcations in $\nu$ around some fixed value $\nu_0$. We therefore fix parameters $(\alpha,\beta,\theta_1,\theta_2,\nu_0)$ and introduce a bifurcation parameter $\mu$ by writing $\nu=\nu_0+\mu$, and we write \eqref{hamvec12}--\eqref{hamvec62} as
\begin{equation}\label{hameqfinal}
\dot{u}=v_{\mathcal{H}^\mu}(u),
\end{equation}
where we have written $\mathcal{H}^\mu$ to indicate that the Hamiltonian depends upon $\mu$.
As mentioned before, $u_0=(0,0,0,0,0,0)$ is an equilibrium solution of \eqref{hameqfinal} and the linearization $L$ of $v_{\mathcal{H}^\mu}(u)$ around this solution, with $\mu=0$, is given by $Lu=(L_1u,L_2u,L_3u,L_4u,L_5u,L_6u)$, where
\begin{align*}
L_1u&=\frac{1}{\beta}\left(\omega+\rho\cos(\theta_1)\int_0^1y\phi_{1y}\ \mathrm{d}y-\cos(\theta_{1})\int_0^1y\phi_{2y}\ \mathrm{d}y\right)-\nu_0\cos(\theta_{1}-\theta_2)\eta_z,\\
L_2u&=\cos(\theta_{1})\int_0^1\left(\psi_1-\rho\cos(\theta_{1})\eta\right)\ \mathrm{d}y-\rho\nu_0\left(\cos(\theta_{1}-\theta_{2})\cos(\theta_{1})-\cos(\theta_{2})\right)\phi_{1z}\vert_{y=1}\\
&\quad -\frac{\cos(\theta_{1})}{h}\int_0^1\left(\psi_2+\cos(\theta_{1})\eta\right)\ \mathrm{d}y+\nu_0\left(\cos(\theta_{1}-\theta_{2})\cos(\theta_{1})-\cos(\theta_{2})\right)\phi_{2z}\vert_{y=1}\\
&\quad -\nu_0^2\beta\sin^2(\theta_{1}-\theta_{2})\eta_{zz}-\nu_0\cos(\theta_{1}-\theta_{2}^0)\omega_{z}+\alpha\eta,\\
L_3u&=\frac{\psi_1}{\rho}-\cos(\theta_{1})\eta-\nu_0\cos(\theta_{1}-\theta_{2})\phi_{1z},\\
L_4u&=-\rho\phi_{1yy}-\nu_0\cos(\theta_{1}-\theta_{2})\psi_{1z}+\frac{\rho\cos(\theta_{1})}{\beta}\left(\omega+\rho\cos(\theta_{1})\int_0^1y\phi_{1y}\ \mathrm{d}y-\cos(\theta_{1})\int_0^1y\phi_{2y}\ \mathrm{d}y\right)\\
&\quad-\rho\nu_0^2\sin^2(\theta_{1}-\theta_{2})\phi_{1zz},\\
L_5u&=\frac{\psi_2}{h}+\frac{\cos(\theta_{1})\eta}{h}-\nu_0\cos(\theta_{1}-\theta_{2})\phi_{2z},\\
L_6u&=-\frac{\phi_{2yy}}{h}-\nu_0\cos(\theta_{1}-\theta_{2})\psi_{2z}-\frac{\cos(\theta_{1})}{\beta}\left(\omega+\rho\cos(\theta_{1})\int_0^1y\phi_{1y}\ \mathrm{d}y-\cos(\theta_{1})\int_0^1y\phi_{2y}\ \mathrm{d}y\right)\\
&\quad-h\nu_0^2\sin^2(\theta_{1}-\theta_{2})\phi_{2zz}.
\end{align*}
where $ \mathcal{D}(L)$ is the set of elements in $X_1$ which satisfy
\begin{align*}
&\begin{cases}
&\phi_{1y}=0,\qquad \qquad \qquad \qquad\qquad \qquad \qquad  \qquad \qquad \qquad  \qquad \qquad y=0,\\
&-\rho\phi_{1y}-\rho\nu_0\left(\cos(\theta_{1}-\theta_{2})\cos(\theta_{1})-\cos(\theta_{2})\right)\eta_z\\
&+\frac{\rho\cos(\theta_{1})}{\beta}\left(\omega+\rho\cos(\theta_{1})\int_0^1y\phi_{1y}\  \mathrm{d}y-\cos(\theta_{1})\int_0^1y\phi_{2y}\  \mathrm{d}y\right)=0,\ y=1,
\end{cases}\\
&\begin{cases}
&\phi_{2y}=0, \qquad \qquad \qquad \qquad\qquad \qquad \qquad  \qquad \qquad \qquad \qquad \qquad y=0,\\
&-\frac{\phi_{2y}}{h}+\nu_0\left(\cos(\theta_{1}-\theta_{2})\cos(\theta_{1})-\cos(\theta_{2})\right)\eta_z\\
&-\frac{\cos(\theta_{1})}{\beta}\left(\omega+\rho\cos(\theta_{1})\int_0^1y\phi_{1y}\  \mathrm{d}y-\cos(\theta_{1})\int_0^1y\phi_{2y}\  \mathrm{d}y\right)=0,\ \ \ y=1.
\end{cases}
\end{align*}
Equation \eqref{hameqfinal} can then be formulated as
\begin{equation}\label{hamiltons_eq}
\dot{u}=Lu+\mathcal{F}^\mu(u),
\end{equation}
where $\mathcal{F}^\mu(u)=v_{\mathcal{H}^\mu}(u)-Lu$.

Finally we note that $v_{\mathcal{H}^\mu}$ anti-commutes with the symmetry
\begin{equation}\label{defreverser}
S:(y,z)\mapsto (y,-z),\ (\eta,\omega,\phi_1,\psi_1,\phi_2,\psi_2)=(\eta,-\omega,-\phi_1,\psi_1,-\phi_2,\psi_2),
\end{equation}
that is, \eqref{hamiltons_eq} is reversible with reverser $S$. The reversibility of \eqref{hamiltons_eq} is due to the invariance of the governing equations \eqref{goveqfinal1}--\eqref{goveqfinal8}, under the transformation 
\begin{equation*}
(x,y,z)\mapsto (-x,y,-z),\ (\eta,\phi_1,\phi_2)\mapsto (\eta,-\phi_1,-\phi_2).
\end{equation*}
Also note that $\mathcal{H}^\mu(Su)=\mathcal{H}^\mu(u)$.
\section{A change of variables}\label{changevar}
The center-manifold theorem applies to equations on linear spaces and so we cannot apply the theorem directly to equation \eqref{hamiltons_eq}, due to the nonlinear boundary conditions \eqref{bdry11}--\eqref{bdry21}.  We therefore make a change of variables in order to obtain an equation equivalent with \eqref{hamiltons_eq}, but with linear boundary conditions. For constructing such variables we follow \cite{Groves2007}.
 
The boundary conditions \eqref{bdry11}--\eqref{bdry21} can be written in the form
\begin{equation}\label{parboundary2}
\phi_{iy}=F_i(u,\mu),\ y=0,1,\ i=0,1,
\end{equation}
where
\begin{align*}
F_1(u,\mu)&=\frac{(1-\eta)y}{\rho}\bigg[-\rho(\nu_0+\mu)^2\sin^2(\theta_1-\theta_2)\eta_z\left(\phi_{1z}-\frac{\eta_z\phi_{1y}}{\eta-1}\right)\\
&\quad-\rho(\nu_0+\mu)\big(\cos(\theta_1-\theta_2)\cos(\theta_1)-\cos(\theta_2)\big)\eta_z\\
&\quad+\sqrt{\frac{(\nu_0+\mu)^2\sin^2(\theta_1-\theta_2)\eta_z^2+1}{\beta^2-\bar{\omega}^2}}\frac{\bar{\omega}(\psi_1-\rho\cos(\theta_1))}{\eta-1}\bigg],\\
F_2(u,\mu)&=(\eta+h)y\bigg[(\nu_0+\mu)^2\sin^2(\theta_1-\theta_2)\eta_z\left(\phi_{2z}-\frac{\eta_z\phi_{2y}}{\eta+h}\right)\\
&\quad+(\nu_0+\mu)\big(\cos(\theta_1-\theta_2)\cos(\theta_1)-\cos(\theta_2)\big)\eta_z\\
&\quad+\sqrt{\frac{(\nu_0+\mu)^2\sin^2(\theta_1-\theta_2)\eta_z^2+1}{\beta^2-\bar{\omega}^2}}\frac{\bar{\omega}(\psi_2-h\cos(\theta_1))}{\eta+h}\bigg].
\end{align*}
Let $V\subseteq  X_1$ be a neighborhood of the origin and let $\Delta$ be a neighborhood of the origin in $\mathbb{R}$. For a fixed value of $\beta$ we choose $V$ small enough so that
\begin{equation*}
-\frac{h}{2}<\eta(z)<\frac{1}{2},\ \abs{\bar{\omega}(z)}<\beta
\end{equation*}
Let $u\in V$, $\mu\in\Delta$ and define $G^{\mu}:V\mapsto X_{1}$, by
\begin{equation*}
G^{\mu}(\eta,\omega,\phi_1,\psi_1,\phi_2,\psi_2)=(\eta,v,\varphi_1,\psi_1,\varphi_2,\psi_2),
\end{equation*}
with 
\begin{align*}
v&=\omega+\int_0^1\rho\cos(\theta_1)y\phi_{1y}\ \mathrm{d}y-\int_0^1\cos(\theta_1)y\phi_{2y}\ \mathrm{d}y,\\
\varphi_1&=\phi_1-\chi_{1y},\\
\varphi_2&=\phi_2-\chi_{2y},
\end{align*}
and where $\chi_i$, $i=1,2$, are the unique solutions of the boundary value problem
\begin{equation*}
\begin{cases}
\chi_{iyy}+\chi_{izz}&=F_i(u,\mu),\\
\chi_i&=0,\qquad\qquad y=0,1.
\end{cases}
\end{equation*}
Note that
\begin{equation*}
\varphi_{iy}=\phi_{iy}-\chi_{iyy}=\phi_{iy}+\chi_{izz}-F_i(u,\mu),\ i=1,2.
\end{equation*}
So if $\phi_i$, $i=1,2$ satisfy \eqref{parboundary2}, then $\varphi_i$, $i=1,2$, satisfy the linear boundary conditions
\begin{equation*}
\varphi_{iy}=0,\ y=0,1,\ i=1,2.
\end{equation*}
The following lemma states that $G^{\mu}$ is a valid change of variables.
\begin{lemma}\label{change_ofvariables_lemma}
\ \begin{itemize}
\item[i]For each $\mu\in\Delta$, the mapping $G^{\mu}$ is a smooth diffeomorphism from the neighborhood $V\subseteq X_1$ of $0$ onto a neighborhood $\tilde{V}\subseteq X_1$ of $0$. The mappings $G^\mu$ and $(G^{\mu})^{-1}$ and their derivatives depend smoothly upon $\mu$.
\item[ii]For each $(u,\mu)\in V\times \Delta$, the operator $\mathrm{d}G^{\mu}[u]:X_1\mapsto X_1$ extends to an isomorphism $\widetilde{\mathrm{d}G}^{\mu}[u]:X_0\mapsto X_0$. The operators $\widetilde{\mathrm{d}G}^{\mu}[u],\ (\widetilde{\mathrm{d}G}^{\mu}[u])^{-1}\in\mathcal{L}(X_0,X_0)$ depend smoothly upon $(u,\mu)\in V\times \Delta$.
\end{itemize}
\end{lemma}
Lemma \ref{change_ofvariables_lemma} can be proven in the same way as \cite[Lemma 3.3]{Groves2007}, by arguing as in \cite[Proposition 2.1]{Bagri2014}. From this change of variables we obtain a Hamiltonian system $(M,\widetilde{\Omega}^{\mu},\widetilde{\mathcal{H}}^\mu)$, where, for $m\in \tilde{V},\ w,w^*\in T_mM,\ \mu\in \Delta$
\begin{equation*}
\widetilde{\Omega}_m^{\mu}(w,w^*)=\Omega\left(\widetilde{\mathrm{d}G}^{\mu}\left[(G^{\mu})^{-1}(m)\right]^{-1}(w),\widetilde{\mathrm{d}G}^{\mu}\left[(G^{\mu})^{-1}(m)\right]^{-1}(w^*)\right),
\end{equation*}
\begin{equation*}
\widetilde{\mathcal{H}}^\mu(m)=\mathcal{H}^\mu((G^{\mu})^{-1}(m)).
\end{equation*}
Hamilton's equation is then given by
\begin{equation}\label{hameqlin}
\dot{w}=v_{\widetilde{\mathcal{H}}^\mu}(w),
\end{equation}
where $v_{\widetilde{\mathcal{H}}^\mu}$ is the Hamiltonian vector field corresponding to the Hamiltonian $\widetilde{\mathcal{H}}^\mu$ and symplectic product $\widetilde{\Omega}^\mu$, with
\begin{equation*}
\mathcal{D}(v_{\widetilde{\mathcal{H}}^\mu})=\{(\eta,v,\varphi_1,\psi_1,\varphi_2,\psi_2)\in \tilde{V}\ : \ \varphi_i=0,\ y=0,1,\ i=0,1\}.
\end{equation*}
Moreover, for elements $w\in\mathcal{D}(v_{\widetilde{\mathcal{H}}^\mu})$ we have that
\begin{equation*}
v_{\widetilde{\mathcal{H}}^{\mu}}(w)=\widetilde{\mathrm{d}G}^{\mu}\left[(G^{\mu})^{-1}(w)\right]\left(v_{\mathcal{H}^\mu}((G^{\mu})^{-1}(w))\right).
\end{equation*}
Let $K$ be the linearization of $v_{\widetilde{\mathcal{H}}^\mu}$ around the equilibrium solution $(0,0,0,0,0,0)$ and $\mu=0$, with
\begin{equation*}
\mathcal{D}(K)=\{(\eta,v,\varphi_1,\psi_1,\varphi_2,\psi_2)\in X_1\ :\ \varphi_i=0,\ y=0,1,\ i=1,2\},
\end{equation*}
so that \eqref{hameqlin} can be written as
\begin{equation}\label{hameqlin2}
\dot{w}=Kw+\widetilde{\mathcal{F}}^\mu(w),
\end{equation}
where $\widetilde{\mathcal{F}}^\mu(w)=v_{\widetilde{\mathcal{H}}^\mu}(w)-Kw$. Note that 
\begin{equation}\label{linopk}
K=\widetilde{\mathrm{d}G}^0[0]L(\mathrm{d}G^0[0])^{-1}.
\end{equation}
Due to this we may work with $L$ instead of $K$ when doing spectral analysis.
\section{Spectrum of $L$}\label{spectrum}
The spectrum of $L$ depends upon the parameters $\alpha,\beta,\theta_1,\theta_2,\nu_0$ and we are interested in parameters for which the number of purely imaginary eigenvalues changes. We will consider two bifurcation scenarios in more detail: a $0^2\i\omega$ resonance and a Hamiltonian-Hopf bifurcation involving mode $\pm 1$ eigenvalues.
 
Let $(\eta,\omega,\phi_1,\psi_1,\phi_2,\psi_2)\in \mathcal{D}(L)$. We expand these functions in Fourier series:
\begin{alignat*}{2}
&\eta(z)=\sum_{k\in\mathbb{Z}}\eta_k\exp(\i kz),&& \omega(z)=\sum_{k\in\mathbb{Z}}\omega_k\exp(\i kz),\\
&\phi_j(y,z)=\sum_{k\in\mathbb{Z}}\phi_{jk}(y)\exp(\i kz),\quad &&\psi_j(y,z)=\sum_{k\in\mathbb{Z}}\psi_{jk}(y)\exp(\i kz),\ j=1,2.
\end{alignat*}
Consider the eigenvalue equation $Lu=\lambda u$. Using the Fourier series expansions above we find that $\lambda\in \mathbb{C}$ is a mode $k$ eigenvalue if and only if
\begin{equation*}
\frac{\rho\left(\nu_0 \i k\cos(\theta_2)+\lambda\cos(\theta_1)\right)^2}{\tan(\gamma_k)}+\frac{\left(\nu_0 \i k\cos(\theta_2)+\lambda \cos(\theta_1)\right)^2}{\tan(h\gamma_k)}=(\alpha-\beta\gamma_k^2)\gamma_k,
\end{equation*}
where
\begin{equation*}
\gamma_k^2=\lambda^2+2\i k\nu_0\lambda\cos(\theta_1-\theta_2)-k^2\nu_0^2.
\end{equation*}
Setting $\lambda=\i s$ we obtain the dispersion relation
\begin{equation}
\frac{\rho\left(\nu_0 k\cos(\theta_2)+s\cos(\theta_1)\right)^2}{\tanh(\tilde{\gamma}_k)}+\frac{\left(\nu_0 k\cos(\theta_2)+s\cos(\theta_1)\right)^2}{\tanh(h\tilde{\gamma}_k)}=(\alpha+\beta\tilde{\gamma}_k^2)\tilde{\gamma}_k,\label{disprel1}
\end{equation}
where
\begin{equation*}
\tilde{\gamma}_k^2=s^2+2k\nu_0 s\cos(\theta_1-\theta_2)+k^2\nu_0^2.
\end{equation*}
We note here that $\i s$ is a mode $k$ eigenvalue if and only if $-\i s$ is a mode  $-k$ eigenvalue. This implies in particular that if $0$ is a mode $k$ eigenvalue then it is also a mode $-k$ eigenvalue. In order to further study higher mode eigenvalues we use the same geometric approach as in \cite{Groves2003} and \cite{HaragusCourcelle2000}.
Let
\begin{align}
l_1&=\nu_0 k\cos(\theta_2)+s\cos(\theta_1),\label{l1}\\
l_2&=\nu_0 k\sin(\theta_2)+s\sin(\theta_1)\label{l2},
\end{align}
and note that $l_1^2+l_2^2=\tilde{\gamma}_k^2$. The dispersion relation \eqref{disprel1} can then be written as
\begin{equation}
\l_1^2\left(\frac{\rho}{\tanh(\sqrt{l_1^2+l_2^2})}+\frac{1}{\tanh(h\sqrt{l_1^2+l_2^2})}\right)-\left(\alpha+\beta(l_1^2+l_2^2)\right)\sqrt{l_1^2+l_2^2}=0.\label{disprel2}
\end{equation}
Denote the left hand side of \eqref{disprel2} by $D(l_1,l_2)$. We see that $\i s$ is a mode $k$ eigenvalue of $L$ if and only if $D(l_1,l_2)=0$ with $l_1$, $l_2$ given by \eqref{l1}, \eqref{l2}. So we are looking for intersections between the set
\begin{equation*}
C_{dr}=\Bigg\{(l_1,l_2)\in \mathbb{R}^2\ :\ l_1^2=\frac{(\alpha+\beta a^2)a}{\frac{\rho}{\tanh(a)}+\frac{1}{\tanh(ha)}},\ l_2^2=a^2-\frac{(\alpha+\beta a^2)a}{\frac{\rho}{\tanh(a)}+\frac{1}{\tanh(ha)}},\ a\in \mathbb{R}\Bigg\}
\end{equation*}
and the lines 
\begin{equation*}
Q_k=\{(l_1,l_2)\ :\ l_1=\nu_0 k\cos(\theta_2)+s\cos(\theta_1),\ l_2=\nu_0 k\sin(\theta_2)+s\sin(\theta_1),\ s\in \mathbb{R}\}.
\end{equation*}
The fact that $\i s$ is a mode $k$ eigenvalue if and only if $-\i s$ is a mode $-k$ eigenvalue can be recovered by noting that $Q_k$ intersect $C_{dr}$ if and only if $Q_{-k}$ intersect $C_{dr}$.

\begin{figure}\centering
\includegraphics{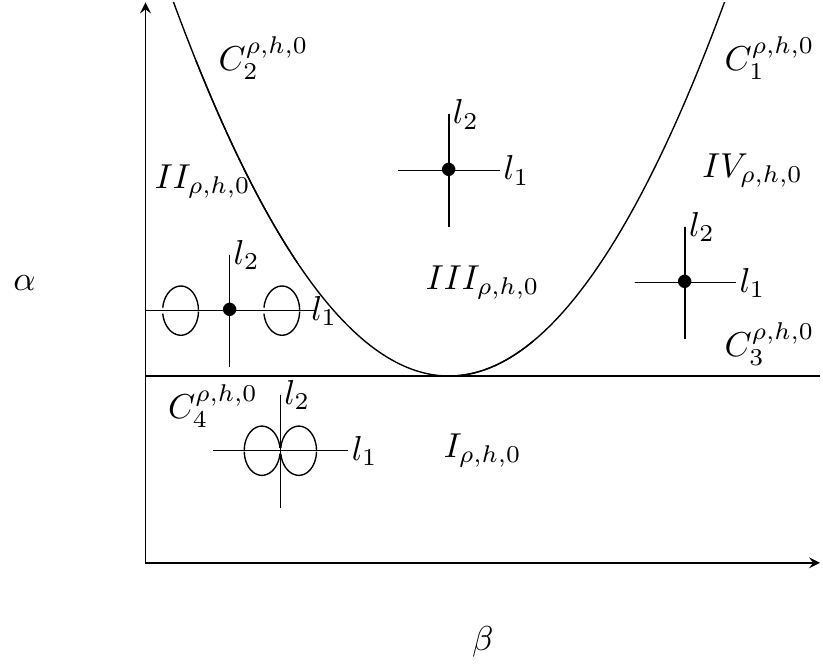}
\caption{Diagram of how the real branches of $C_{dr}$ look in certain regions of the $(\beta,\alpha)$-plane.}
\label{bifdiagrambranches}
\end{figure}

In Figure \ref{bifdiagrambranches} we provide a simplified picture $C_{dr}$ for different values of $\alpha$ and $\beta$. Qualitatively, the same picture was also obtained for surface waves in \cite{HaragusCourcelle2000}. Note in particular that for $(\beta,\alpha)\in III_{\rho,h,0}\cup IV_{\rho,h,0}$ there are no purely imaginary nontrivial mode $k$ eigenvalues. The curves $C_i^{\rho,h,\theta_1}$, $i=1,..,4$ are defined in \eqref{c1}--\eqref{c4}. Note that the curves $C_i^{\rho,h,0}$ obtained when choosing $\theta_1=0$ are the bifurcation curves found in the study of two-dimensional internal solitary waves (see \cite{nilsson2016}).
In Figure \ref{inter2} we have sketched the case when $Q_0$ intersects $C_{dr}$ at the origin and in two additional points and $Q_{\pm 1}$ intersect $C_{dr}$ in two distinct points each.
We see from this picture that in general there are at least two different mode $\pm 1$ eigenvalues or there are none. For the special case when $Q_k$ is tangent to $C_{dr}$, $\i s$ is a mode $k$ eigenvalue of algebraic multiplicity $2$. More precisely, we find from equation \eqref{disprel1} that $\i s$, with $s+kv_0\cos(\theta_1-\theta_2)\neq 0$, is a mode $k$ eigenvalue of algebraic multiplicity 2 if and only if
\begin{align}
\alpha=\alpha_k^*(s)&:=-\beta_k^*(s)\tilde{\gamma}_k^2+\frac{1}{\tilde{\gamma}_k}\left(s\cos(\theta_1)+\nu_0 k\cos(\theta_2)\right)^2\left(\frac{\rho}{\tanh(\tilde{\gamma}_k)}+\frac{1}{\tanh(h\tilde{\gamma}_k)}\right),\label{parametrisationalpha}\\
\beta=\beta_k^*(s)&:=\frac{s\cos^2(\theta_1)+k\nu_0\cos(\theta_1)\cos(\theta_2)}{\tilde{\gamma}_k(s+\nu_0 k\cos(\theta_1-\theta_2))}\left(\frac{\rho}{\tanh(\tilde{\gamma}_k)}+\frac{1}{\tanh(h\tilde{\gamma}_k)}\right)\nonumber\\
&\quad -\frac{(s\cos(\theta_1)+\nu_0 k\cos(\theta_2))^2}{2\tilde{\gamma}_k^3}\left[\frac{\rho}{\tanh(\tilde{\gamma}_k)}+\frac{1}{\tanh(h\tilde{\gamma}_k)}\right.\nonumber\\
&\qquad\qquad\qquad\qquad\qquad\qquad\quad\ \ \left.+\tilde{\gamma}_k\left(\frac{\rho}{\sinh^2(\tilde{\gamma}_k)}+\frac{h}{\sinh^2(h\tilde{\gamma}_k)}\right)\right] \label{parametrisationbeta}.
\end{align}
In the special case when $s=-kv_0\cos(\theta_1-\theta_2)$, the dispersion relation \eqref{disprel1} reduces to
\begin{align}
(\alpha+\beta \nu_0^2k^2\sin^2(\theta_1-\theta_2))\nu_0k\sin(\theta_1-\theta_2)&= \nu_0^2k^2\sin^2(\theta_1)\sin^2(\theta_1-\theta_2)\bigg(\frac{\rho}{\tanh(\nu_0k\sin(\theta_1-\theta_2))}\nonumber\\
&\quad\quad\quad+\frac{1}{\tanh(h\nu_0k\sin(\theta_1-\theta_2))}\bigg),\label{0eigenmodek}
\end{align}
and we find that $\i s$ is of algebraic multiplicity $2$ if and only if
\begin{equation}\label{condition0}
\nu_0k\cos(\theta_1)\sin(\theta_1)=0.
\end{equation}
We note that \eqref{condition0} is satisfied for $k\neq 0$ if and only if $\theta_1=0$ or $\theta_1=\pm\pi/2$. However, equation \eqref{0eigenmodek} has no solutions when $k\neq 0$ and $\theta_1=0$. Hence, $\i s=-\i \nu_0k\cos(\theta_1-\theta_2)$ is an eigenvalue of algebraic multiplicity $2$ if and only if $\theta_1=\pm \pi/2$ and $(\beta,\alpha)$ belongs to the line defined by \eqref{0eigenmodek}. In addition, when $\theta_1=\pm \pi/2$ then an algebraically double mode $k$ eigenvalue $\i s$ must necessarily satisfy $s=-\nu_0 k\cos(\pm \pi/2-\theta_2)$. Indeed, when $\theta_1=\pm \pi/2$, then $Q_k$ is parallel with the $l_2$-axis, so $Q_k$ is tangent to $C_{dr}$ only when $l_2=0$, which implies that $s=\pm\nu_0k\sin(\theta_2)=-\nu_0k\cos(\pm\pi/2-\theta_2)$. The case when $k=0$ will be investigated more thoroughly in section \ref{mode0}.

In the further special case $(\theta_1,\theta_2)=(\pm \pi/2,0)$ it follows from the dispersion relation \eqref{disprel1} that $\i s$ is a mode $k$ eigenvalue if and only if it is a mode $-k$ eigenvalue. Hence, all eigenvalues are of geometric multiplicity $2$. Also note that since $\theta_1=\pi/2$, we have that an algebraically double mode $k$ eigenvalue $\i s$ must satisfy $s=-\nu_0k\cos(0\pm\pi/2)=0$.
Similarly, in the other special case when $(\theta_1,\theta_2)=(0,\pm\pi/2)$ we again have that all mode $k$ eigenvalues are of geometric multiplicity $2$. Both of these cases have been studied in the surface wave setting. The case when $(\theta_1,\theta_2)=(0,\pm\pi/2)$ was considered in \cite{Groves2007} and the case when $(\theta_1,\theta_2)=(\pm\pi/2,0)$ was considered in \cite{Groves2001,HK2001}. Moreover, both of these cases were again investigated in \cite{Groves2003}.

\begin{figure}\centering
\includegraphics{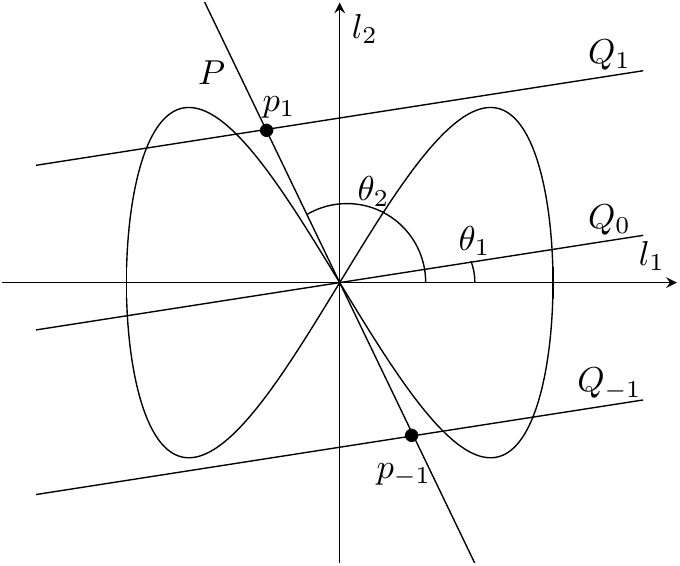}
\caption{Intersections of P and $Q_{\pm 1}$.}
\label{LQ0}
\end{figure}

The following characterization of the purely imaginary eigenvalues, used also in \cite{Groves2003}, will be helpful when discussing the $00(\i s)(\i \kappa_0)$ resonance in section \ref{mode1}. Note that
\begin{equation*}
Q_k=\{sQ_0+\nu_0 kP\ :\ s\in \mathbb{R}\},
\end{equation*}
where $P=(\cos(\theta_2),\sin(\theta_2))$ and $Q_0$ is here to be interpreted as $(\cos(\theta_1),\sin(\theta_1))$. This means that points on the lines $Q_k$ have coordinates $(s,\nu_0k)$ in the coordinate system $(Q_0,P)$. The line generated by $P$ intersects the lines $Q_k$ in points $p_k=(\nu_0k\cos(\theta_2),\nu_0k\sin(\theta_2))$, see Figure \ref{LQ0}. The imaginary part of a purely imaginary mode $k$ eigenvalue can therefore be interpreted as the signed distance between $p_k$ and the corresponding intersection of $Q_k$ and $C_{dr}$. 

\subsection{Mode $0$ eigenvalues}\label{mode0}
\begin{figure}\centering
\includegraphics{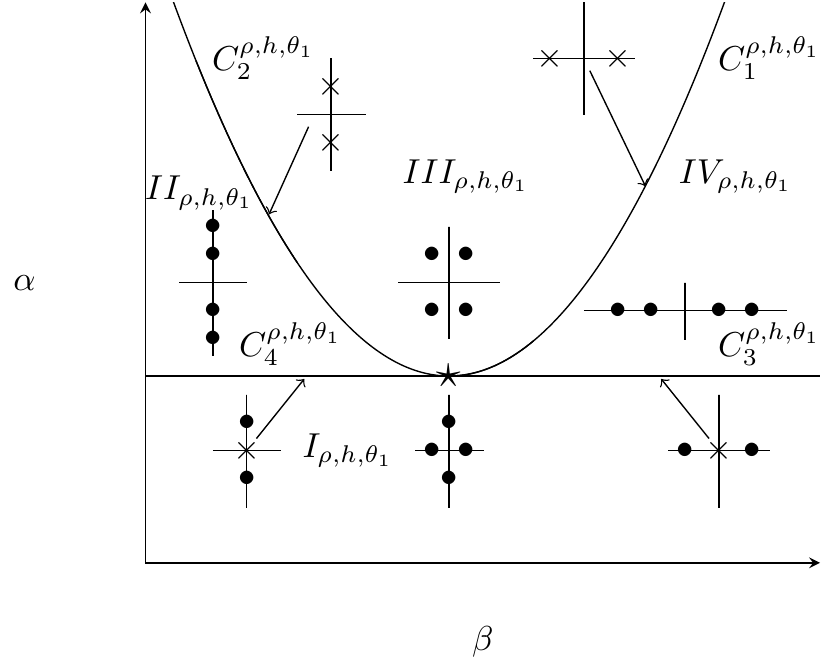}
\caption{Bifurcation curves in the $(\beta,\alpha)$-plane for mode $0$ eigenvalues. The $\star$ in the picture indicates the point $\cos^2(\theta_1)(\frac{\rho+h}{3},\rho+\frac{1}{h})$.}
\label{fig2}
\end{figure}

For $k=0$ \eqref{disprel1} becomes
\begin{equation}\label{disprelmode0}
s^2\cos^2(\theta_1)\left(\frac{\rho}{\tanh(s)}+\frac{1}{\tanh(hs)}\right)=(\alpha+\beta s^2)s.
\end{equation}
Note that $0$ is a solution of \eqref{disprelmode0}. In fact $0$ has for all parameter values, two eigenvectors, each with a corresponding generalized eigenvector, that is, $0$ is trivially a mode $0$ eigenvalue of algebraic multiplicity $4$. Also note that when $\theta_1=\pm\pi/2$ there are no other purely imaginary mode $0$ eigenvalues.
From \eqref{disprelmode0} we obtain the bifurcation diagram in Figure \ref{fig2}.
The curves in Figure \ref{fig2} are
\begin{align}
C_1^{\rho,h,\theta_1}&=\left\{(\beta,\alpha)=(\beta^*(\i s),\alpha^*(\i s)),\ s\in(0,\min(0,\pi/h))\right\},\label{c1}\\
C_2^{\rho,h,\theta_1}&=\left\{(\beta,\alpha)=(\beta^*(s),\alpha^*(s)),\ s\in(0,\infty)\right\},\label{c2}\\
C_3^{\rho,h,\theta_1}&=\left\{(\beta,\alpha)\ :\ \beta>\cos^2(\theta_1)\frac{\rho+h}{3},\ \alpha=\cos^2(\theta_1)\left(\rho+\frac{1}{h}\right)\right\}\label{c3},\\
C_4^{\rho,h,\theta_1}&=\left\{(\beta,\alpha)\ :\ \beta<\cos^2(\theta_1)\frac{\rho+h}{3},\ \alpha=\cos^2(\theta_1)\left(\rho+\frac{1}{h}\right)\right\}.\label{c4}
\end{align}
\begin{itemize}
\item When the curve $C_1^{\rho,h,\theta_1}$ is crossed from below, the mode $0$ eigenvalues of $L$ change from two pairs of real simple eigenvalues to a plus-minus complex-conjugate quartet of complex simple eigenvalues. For points on the curve $C_1^{\rho,h,\theta_1}$ the eigenvalues collide on the real axis to form a plus-minus pair of algebraically double real eigenvalues. This change in the spectrum is called a Hamiltonian real 1:1 resonance.
\item When the curve $C_2^{\rho,h,\theta_1}$ is crossed from below, the mode $0$ eigenvalues of $L$ change from two pairs of simple imaginary eigenvalues to a plus-minus complex-conjugate quartet of complex simple eigenvalues. For points on the curve $C_2^{\rho,h,\theta_1}$ the eigenvalues collide on the imaginary axis to form a plus-minus pair of algebraically double imaginary eigenvalues. This change in the spectrum is called a Hamiltonian-Hopf bifurcation.
\item When the curve $C_3^{\rho,h,\theta_1}$ is crossed from below, the mode $0$ eigenvalues of $L$ change from a pair of simple real eigenvalues and a pair of simple imaginary eigenvalues to two pairs of simple real eigenvalues. For points on the curve $C_3^{\rho,h,\theta_1}$ the imaginary eigenvalues collide at $0$, making $0$ an eigenvalue of algebraic multiplicity $6$. This change in the spectrum is called a Hamiltonian real $0^2$ resonance. 
\item When the curve $C_4^{\rho,h,\theta_1}$ is crossed from below, the mode $0$ eigenvalues of $L$ change from a pair of simple real eigenvalues and a pair of simple imaginary eigenvalues to two pairs of simple imaginary eigenvalues. For points on the curve $C_4^{\rho,h.\theta_1}$ real eigenvalues collide at $0$, making $0$ an eigenvalue of algebraic multiplicity $6$.
\item When $(\beta,\alpha)=\cos^2(\theta_1)(\frac{\rho+h}{3},\rho+\frac{1}{h})$, $0$ is an eigenvalue of algebraic multiplicity $8$, and there are no other purely imaginary mode $0$ eigenvalues.
\end{itemize}
\subsection{Bifurcations involving mode $\pm 1$ eigenvalues}\label{mode1}
We begin by describing the $00(\i s)(\i \kappa_0)$ resonance mentioned in the introduction. Here we want to choose $\nu_0$ so that $\pm \i \kappa_0,0$ are contained in the spectrum of $L$, where $\pm \i \kappa_0$ are non-zero, simple mode $0$ eigenvalues, $0$ is a simple mode $\pm 1$ eigenvalue (and trivially a mode $0$ eigenvalue of algebraic multiplicity 4). We allow for other purely imaginary mode $\pm k$ eigenvalues $\pm \i s$ under the nonresonance condition $s\neq m\kappa_0$ for all $m\in \mathbb{Z}$. In order for $\kappa_0\neq 0$, we need to assume that $\theta_1\neq \pm \pi/2$. The mode $0$ eigenvalues of $L$ were investigated in section \ref{mode0}; we therefore turn to the mode $\pm 1$ eigenvalue $0$.
From \eqref{disprel1} we get that $0$ is a mode $\pm 1$ eigenvalue if and only if
\begin{figure}\centering
\includegraphics{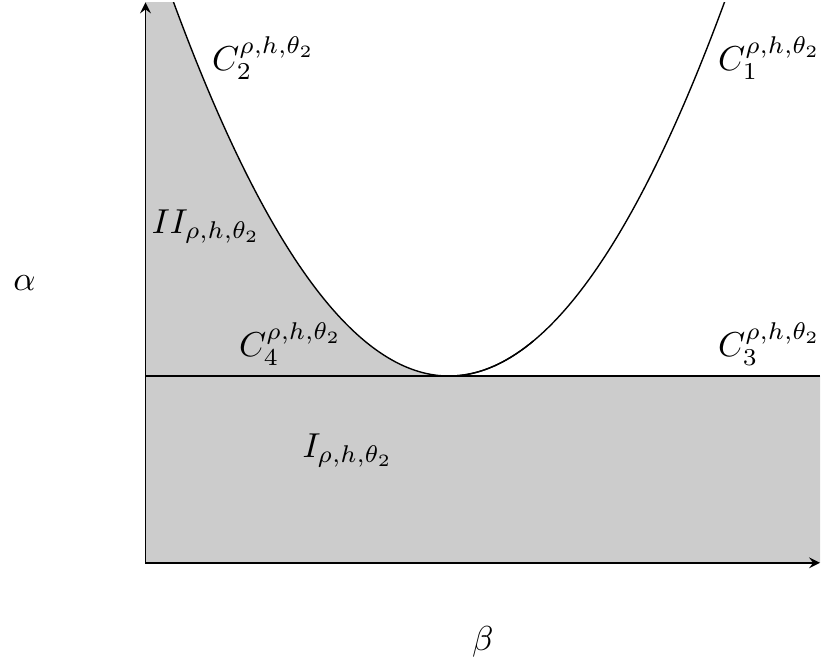}
\caption{Bifurcation curves related to equation \eqref{disprel0mode1}. In region $I_{\rho,h,\theta_2}$ there is one positive solution $\nu_0$ of \eqref{disprel0mode1} and in $II_{\rho,h,\theta_2}$ there are two such solutions.}
\label{fig3}
\end{figure}
\begin{equation}\label{disprel0mode1}
\nu_0^2\cos^2(\theta_2)\left(\frac{\rho}{\tanh(\nu_0)}+\frac{1}{\tanh(h\nu_0)}\right)=(\alpha+\beta \nu_0^2)\nu_0.
\end{equation}
Here we need to assume that $\theta_2\neq\pm \pi/2$ or else there are no nontrivial solutions $\nu_0$ of \eqref{disprel0mode1}. We note that \eqref{disprel0mode1} is analogous to \eqref{disprelmode0} and we can therefore conclude that it possesses real solutions if $(\beta,\alpha)$ belong to the shaded region in Figure \ref{fig3}. More specifically, for $(\beta,\alpha)\in I_{\rho,h,\theta_2}$ there exists one positive solution $\nu_0$ of \eqref{disprel0mode1}, and for $(\beta,\alpha)\in II_{\rho,h,\theta_2}$ there exist two positive solutions $\nu_1,\nu_2$. There are no solutions of \eqref{disprel0mode1} for values of $(\beta,\alpha)$ belonging to any of  the remaining regions of the $(\beta,\alpha)$-plane. So given $(\beta,\alpha)\in I_{\rho,h,\theta_2}$ we obtain a positive solution $\nu_0$ of \eqref{disprel0mode1} and for this choice of $\nu_0$, $0$ is a mode $\pm 1$ eigenvalue of $L$. Given $(\beta,\alpha)\in II_{\rho,h,\theta_2}$ we obtain two positive solutions of \eqref{disprel0mode1}, and we choose one of these solutions as our $\nu_0$. The relevant transition curves in Figure \ref{fig3} are
\begin{align*}
C_2^{\rho,h,\theta_2}&=\left\{(\beta,\alpha)=(\tilde{\beta}(\nu_0),\tilde{\alpha}(\nu_0)\right\},\\
C_4^{\rho,h,\theta_2}&=\left\{(\beta,\alpha)\ :\ \beta<\cos^2(\theta_2)\left(\frac{\rho+h}{3}\right),\ \alpha=\cos^2(\theta_2)\left(\rho+\frac{1}{h}\right)\right\},
\end{align*}
where
\begin{align*}
\tilde{\alpha}(\nu_0)&=-\tilde{\beta}(\nu_0)\nu_0^2+\nu_0\cos^2(\theta_2)\left(\frac{\rho}{\tanh(\nu_0)}+\frac{1}{\tanh(h\nu_0)}\right),\\
\tilde{\beta}(\nu_0)&=\frac{\cos^2(\theta_2)}{2\nu_0}\left(\frac{\rho}{\tanh(\nu_0)}+\frac{1}{\tanh(h\nu_0)}\right)-\frac{\cos^2(\theta_2)}{2}\left(\frac{\rho}{\sinh^2(\nu_0)}+\frac{1}{\sinh^2(h\nu_0)}\right).
\end{align*}
These are essentially the same curves as $C_2^{\rho,h,\theta_1}$, $C_4^{\rho,h,\theta_1}$, in fact
\begin{equation*}
(\tilde{\beta}(\nu_0),\tilde{\alpha}(\nu_0))=\frac{\cos^2(\theta_2)}{\cos^2(\theta_1)}(\beta^*(\nu_0),\alpha^*(\nu_0)).
\end{equation*}
Note that we can allow values $(\beta,\alpha)$ belonging to either $C_2^{\rho,h,\theta_2}$ or $C_4^{\rho,h,\theta_2}$ since this does not change the multiplicity of the eigenvalue $0$. In fact, under the assumption that $\theta_1\neq\pm\pi/2$, we know from section \ref{spectrum} that $0$ is a mode $\pm 1$ eigenvalue of algebraic multiplicity $2$ if and only if $(\beta,\alpha)=(\beta_1^*(0),\alpha_1^*(0))$.  In conclusion we find that the following holds.
If $\cos^2(\theta_1)\leq \cos^2(\theta_2)$, then $I_{\rho,h,\theta_1}\subseteq I_{\rho,h,\theta_2}$. If $(\beta,\alpha)\in I_{\rho,h,\theta_1}$, then we know from section \ref{mode0} that $L$ has the mode $0$ eigenvalues $\pm \i\kappa_0$ and no other nontrivial mode $0$ eigenvalues. Since $(\beta,\alpha)\in I_{\rho,h,\theta_2}$ as well, there exists a solution $\nu_0$ of \eqref{disprel0mode1}, which implies that $0$ is a mode $\pm 1$ eigenvalue of $L$, for this choice of $\nu_0$. If $(\beta,\alpha)\in II_{\rho,h,\theta_1}$, we get from section \ref{mode0} that $L$ has the mode $0$ eigenvalues $\pm \i\kappa_0$, $\pm \i\kappa_1$, which are generically nonresonant. In addition $(\beta,\alpha)$ belongs to either $I_{\rho,h,\theta_2}$ or $II_{\rho,h,\theta_2}$. In the first case, equation \eqref{disprel0mode1} has, as previously mentioned, one solution $\nu_0$ and in the second case there are two solutions $\nu_1,\nu_2$. Hence, in both cases we can choose $\nu_0$ so that $0$ is a mode $\pm 1$ eigenvalue of $L$, where we in the second case choose either $\nu_1$ or $\nu_2$ and use this as our $\nu_0$. If $\cos^2(\theta_2)\leq \cos^2(\theta_1)$ we can argue in the same way, but with $\theta_1$ and $\theta_2$ reversed.
\begin{figure}\centering
\includegraphics{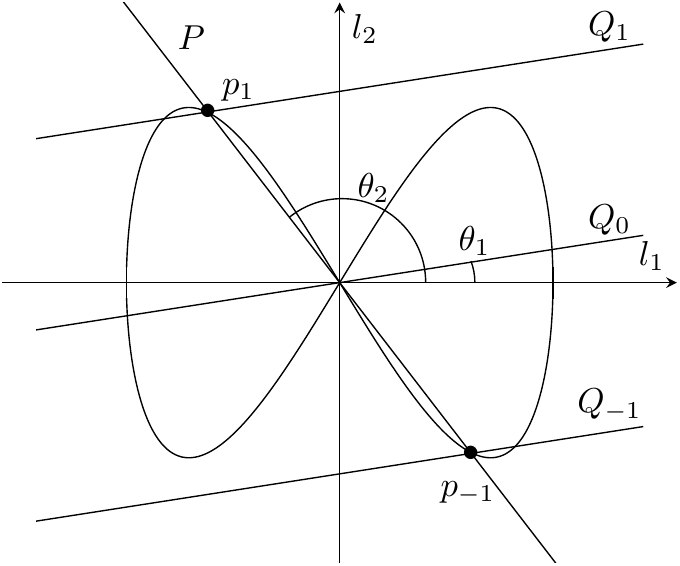}
\caption{Intersections between $Q_{\pm 1}$ when $0$ is a mode $\pm 1$ eigenvalue.} 
\label{0eigenmodepm1}
\end{figure}
This bifurcation can also be explained using the geometric interpretation of the dispersion relation.  First note that the angle between $Q_0$ and the positive $l_1$ axis is $\theta_1$ and the angle between $P$ and the positive $l_1$ axis is $\theta_2$. In order for $0$ to be a mode $\pm 1$ eigenvalue we need to have that the points $p_{\pm 1}$ and the intersection points of $Q_{\pm 1}$ and $C_{dr}$ are the same, see Figure \ref{0eigenmodepm1}. In particular we need to ensure that $P$ intersects $C_{dr}$. In Figure \ref{0eigenmodepm1}, such an intersection is possible if the angle that $C_{dr}$ makes with the negative $l_1$ axis at the origin, in the second quadrant, is greater than $\pi-\theta_2$. In region $I_{\rho,h,0}$ the angle between $C_{dr}$ and the negative $l_1$ axis at the origin, is given by $\pm \arctan((\frac{\rho+1/h}{\alpha}-1)^{\frac{1}{2}})$. Hence, by choosing $\alpha$ small enough we can ensure that $P$ intersects $C_{dr}$. We also see from Figure \ref{0eigenmodepm1} that $L$ will have at least one other pair of simple mode $\pm 1$ eigenvalues $\pm\i s$ in addition to the mode $\pm 1$ eigenvalue $0$.

\begin{figure}\centering
\includegraphics{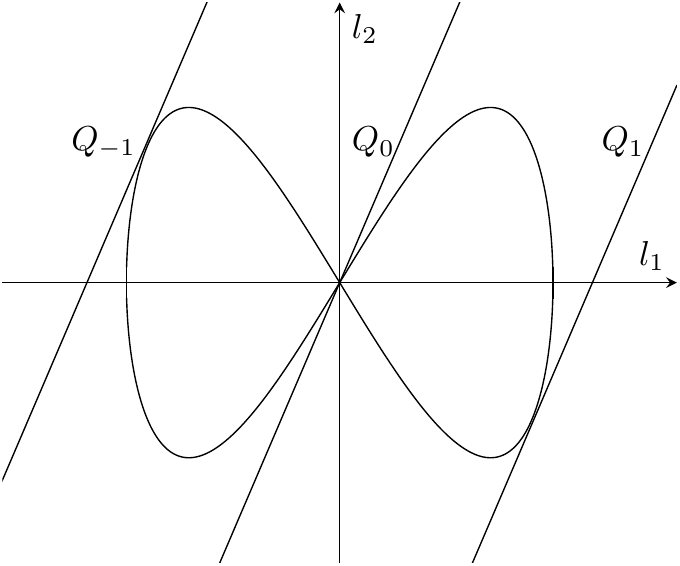}
\caption{Intersections between the real branches of $C_{dr}$ and the lines $Q_{\pm 1}$, $Q_0$ for a critical value $\nu_0$.}
\label{inter1}
\end{figure}

We next consider the case when a Hamiltonian-Hopf bifurcation occurs, involving mode $\pm 1$ eigenvalues. If $\nu_0$ is chosen sufficiently large, the only line $Q_k$ which can intersect $C_{dr}$ is $Q_0$, and so there are no higher mode eigenvalues of $L$. If $\nu_0$ is then decreased there will be some critical value of $\nu_0$ such that the lines $Q_1$, $Q_{-1}$ are tangent to $C_{dr}$, which means that $L$ has the mode $\pm 1$ eigenvalues $\pm \i s$ of algebraic multiplicity 2. This case is illustrated in Figure \ref{inter1}. If $\nu_0$ is decreased further, then we have the case illustrated in Figure \ref{inter2}, that is $L$ has two simple mode $\pm 1$ eigenvalues. This shows that a Hamiltonian-Hopf bifurcation occurs at some critical value of $\nu_0$. We will focus on the case when $L$ does not have any other nontrivial imaginary eigenvalues. This is achieved for values of $(\beta,\alpha)$ belonging to $III_{\rho,h,\theta_1}$ or $IV_{\rho,h,\theta_1}$. Recall from  Figure \ref{bifdiagrambranches} that for $(\beta,\alpha)\in III_{\rho,h,0}\cup IV_{\rho,h,0}$ $L$ has no purely imaginary eigenvalues. We must therefore assume that $\theta_1\neq 0$.
\subsection{Generalized eigenvectors}\label{geneigen}
The eigenvector corresponding to the mode $k$ eigenvalue $\i s$, is given by $\exp(\i kz)v_s^k$, with
\begin{equation*}
v_{s}^k=\left(\begin{array}{c}
\frac{\tilde{\gamma}_{k}}{k\nu_0\cos(\theta_2)+s\cos(\theta_1)}\\
-\i \rho\cos(\theta_1)\left(\frac{1}{\tanh(\tilde{\gamma}_{k})}-\frac{1}{\tilde{\gamma}_{k}}\right)-\i\cos(\theta_1)\left(\frac{1}{\tanh(h\tilde{\gamma}_{k})}-\frac{1}{h\tilde{\gamma}_{k}}\right)+\frac{\i\tilde{\gamma}_{k}\beta(k\nu_0\cos(\theta_1-\theta_2)+s)}{k\nu_0\cos(\theta_2)+s\cos(\theta_1)}\\
\frac{\i\cosh(\tilde{\gamma}_{k}y)}{\sinh(\tilde{\gamma}_{k})}\\
\rho\left(\frac{\tilde{\gamma}_{k}\cos(\theta_1)}{k\nu_0\cos(\theta_2)+s\cos(\theta_1)}-\frac{\cosh(\tilde{\gamma}_{k}y)}{\sinh(\tilde{\gamma}_{k})}(k\nu_0\cos(\theta_1-\theta_2)+s)\right)\\
-\frac{\i\cosh(h\tilde{\gamma}_{k}y)}{\sinh(h\tilde{\gamma}_{k})}\\
-\frac{\tilde{\gamma}_{k}\cos(\theta_1)}{k\nu_0\cos(\theta_2)+s\cos(\theta_1)}+\frac{\cosh(h\tilde{\gamma}_{k}y)}{\sinh(h\tilde{\gamma}_{k})}h(k\nu_0\cos(\theta_1-\theta_2)+s)
\end{array}\right).
\end{equation*}
If $\i s$ is a mode $k$ eigenvalue of algebraic multiplicity $2$ there is a generalized eigenvector $\exp(\i kz)u_s^k$, where $u_s^k=(\eta_k,\omega_k,\phi_{1k},\psi_{1k},\phi_{2k},\psi_{2k})$, such that $(L-\i s)\exp(\i kz)u_s^k=\exp(\i kz)v_s^k$, where
\begin{align*}
\eta_k&=0,\\
\omega_k&=\frac{\rho\cos(\theta_1)(s+k\nu_0\cos(\theta_1-\theta_2))}{\tilde{\gamma}_k^2}\left(\frac{1}{\tanh(\tilde{\gamma}_k)}+\frac{\tilde{\gamma}_k}{\sinh^2(\tilde{\gamma}_k)}-\frac{2}{\tilde{\gamma}_k}\right)\\
&\quad+\frac{\rho\cos^2(\theta_1)}{k\nu_0\cos(\theta_2)+s\cos(\theta_1)}\left(\frac{1}{\tilde{\gamma}_k}-\frac{1}{\tanh(\tilde{\gamma}_k)}\right)\\
&\quad+\frac{\cos(\theta_1)(s+k\nu_0\cos(\theta_1-\theta_2))}{\tilde{\gamma}_k^2}\left(\frac{1}{\tanh(h\tilde{\gamma}_k)}+\frac{h\tilde{\gamma}_k}{\sinh^2(h\tilde{\gamma}_k)}-\frac{2}{h\tilde{\gamma}_k}\right)\\
&\quad +\frac{\cos^2(\theta_1)}{k\nu_0\cos(\theta_2)+s\cos(\theta_1)}\left(\frac{1}{h\tilde{\gamma}_k}-\frac{1}{\tanh(h\tilde{\gamma}_k)}\right)+\frac{\tilde{\gamma}_k\beta}{\nu_0k\cos(\theta_2)+s\cos(\theta_1)},\\
\phi_{1k}&=\frac{s+k\nu_0\cos(\theta_1-\theta_2)}{\tilde{\gamma}_k^2\sinh(\tilde{\gamma}_k)}\left(\tilde{\gamma}_ky\sinh(\tilde{\gamma}_ky)-\left(1+\frac{\tilde{\gamma}_k}{\tanh(\tilde{\gamma}_k)}\right)\cosh(\tilde{\gamma}_ky)\right)\\
&\quad+\frac{\cos(\theta_1)\cosh(\tilde{\gamma}_ky)}{(k\nu_0\cos(\theta_2)+s\cos(\theta_1))\sinh(\tilde{\gamma}_k)},\\
\psi_{1k}&=\i\rho(s+k\nu_0\cos(\theta_1-\theta_2))\bigg[\frac{s+k\nu_0\cos(\theta_1-\theta_2)}{\tilde{\gamma}_k^2\sinh(\tilde{\gamma}_k)}\left(\tilde{\gamma}_ky\sinh(\tilde{\gamma}_ky)-\left(1+\frac{\tilde{\gamma}_k}{\tanh(\tilde{\gamma}_k)}\right)\cosh(\tilde{\gamma}_ky)\right)\\
&\quad+\frac{\cos(\theta_1)\cosh(\tilde{\gamma}_ky)}{(k\nu_0\cos(\theta_2)+s\cos(\theta_1))\sinh(\tilde{\gamma}_k)}\bigg]+\frac{\i\rho\cosh(\tilde{\gamma}_ky)}{\sinh(\tilde{\gamma}_k)},\\
\phi_{2k}&=-\frac{s+k\nu_0\cos(\theta_1-\theta_2)}{\tilde{\gamma}_k^2\sinh(h\tilde{\gamma}_k)}\left(h\tilde{\gamma}_ky\sinh(h\tilde{\gamma}_k)-\left(1+\frac{h\tilde{\gamma}_k}{\tanh(h\tilde{\gamma}_k)}\right)\cosh(h\tilde{\gamma}_k)\right)\\
&\quad-\frac{\cos(\theta_1)\cosh(h\tilde{\gamma}_ky)}{(k\nu_0\cos(\theta_2)+s\cos(\theta_1))\sinh(h\tilde{\gamma}_k)},\\
\psi_{2k}&=-\i h(s+k\nu_0\cos(\theta_1-\theta_2))\bigg[\frac{s+k\nu_0\cos(\theta_1-\theta_2)}{\tilde{\gamma}_k^2\sinh(h\tilde{\gamma}_k)}\left(h\tilde{\gamma}_ky\sinh(h\tilde{\gamma}_k)-\left(1+\frac{h\tilde{\gamma}_k}{\tanh(h\tilde{\gamma}_k)}\right)\cosh(h\tilde{\gamma}_k)\right)\\
&\quad+\frac{\cos(\theta_1)\cosh(h\tilde{\gamma}_ky)}{(k\nu_0\cos(\theta_2)+s\cos(\theta_1))\sinh(h\tilde{\gamma}_k)}\bigg]-\frac{\i h\cosh(h\tilde{\gamma}_ky)}{\sinh(h\tilde{\gamma}_k)}.
\end{align*}
In addition, for $\theta_1\neq \pm\pi/2$, we find that the mode $k$ eigenvalue $\i s$ is of algebraic multiplicity $3$ if
\begin{equation}\label{highermult}
\frac{\mathrm{d}\beta_k^*(s)}{\mathrm{d}s}=0,
\end{equation}
However, this case will not be investigated further in the present paper.
The mode $0$ eigenvalue $0$ is trivially of algebraic multiplicity $4$, with eigenvectors $e_1,\ e_2$ and corresponding generalized eigenvectors $f_1,\ f_2$, satisfying $Lf_1=e_1,\ Lf_2=e_2$. These are given by
\begin{equation*}
e_1=\left(\begin{array}{c}
0\\
0\\
1\\
0\\
0\\
0
\end{array}\right),\
e_2=\left(\begin{array}{c}
0\\
0\\
0\\
0\\
1\\
0
\end{array}\right),\ 
f_1=\left(\begin{array}{c}
-\frac{\rho\cos(\theta_1)}{\alpha}\\
0\\
0\\
\rho\left(1-\frac{\rho\cos^2(\theta_1)}{\alpha}\right)\\
0\\
\frac{\rho\cos^2(\theta_1)}{\alpha}
\end{array}\right),\ 
f_2=\left(\begin{array}{c}
\frac{\cos(\theta_1)}{\alpha}\\
0\\
0\\
\frac{\rho\cos^2(\theta_1)}{\alpha}\\
0\\
h-\frac{\cos^2(\theta_1)}{\alpha}
\end{array}\right).
\end{equation*}
\section{Center manifold reduction}\label{cmtsection}
We will use the following version of the center-manifold theorem which is due to Mielke \cite{MI} and was used in for example \cite{BG}.
\begin{theorem}\label{cmt}
Consider the differential equation
\begin{equation}
\dot{u}=Ku+\mathcal{F}(u,\mu),\label{fulleqthm}
\end{equation}
where $u$ belongs to a Hilbert space $E$, $\mu\in \mathbb{R}^n$ is a parameter and $K:\mathcal{D}(K)\subset E\mapsto E$ is a closed linear operator. Suppose that \eqref{fulleqthm} is Hamilton's equation for the Hamiltonian system $(E,\Omega,\mathcal{H})$. Suppose further that
\begin{itemize}
\item[H1.] $E$ has two closed, $K$-invariant subspaces $E_1$, $E_2$ such that
\begin{align*}
E&=E_1\oplus E_2,\\
\dot{u}_1&=K_1u_1+\mathcal{F}_1(u_1+u_2,\mu),\\
\dot{u}_2&=K_2u_2+\mathcal{F}_2(u_1+u_2,\mu),
\end{align*}
where $K_i=K\vert_{\mathcal{D}(K)\cap E_i}:\mathcal{D}(K)\cap E_i\mapsto E_i$, $i=1,2$ and $\mathcal{F}_1=P\mathcal{F}$, $\mathcal{F}_2=(I-P)\mathcal{F}$, where $P$ is the projection of $E$ onto $E_1$.
\item[H2.] $E_1$ is finite dimensional and the spectrum of $K_1$ lies on the imaginary axis.
\item[H3.] The imaginary axis lies in the resolvent set of $K_2$ and 
\begin{equation*}
\norm{(K_2-\i aI)^{-1}}_{E_2\mapsto E_2}\leq \frac{C}{1+\abs{a}},\quad a\in \mathbb{R}.
\end{equation*}
\item[H4.]
There exists $k\in \mathbb{N}$ and neighborhoods $\Lambda\subset \mathbb{R}^n$ and $U\subset \mathcal{D}(K)$ of $0$ such that $\mathcal{F}$ is $k+1$ times continuously differentiable on $U\times \Lambda$ and the derivatives of $\mathcal{F}$ are bounded and uniformly continuous on $U\times \Lambda$ with
\begin{equation*}
\mathcal{F}(0,0)=0,\quad d_1\mathcal{F}[0,0]=0.
\end{equation*}
\end{itemize}
Under the hypothesis $H1-H4$ there exist neighborhoods $\tilde{\Lambda}\subset \Lambda$ and $\tilde{U}_1\subset U\cap E_1$, $\tilde{U}_2\subset U\cap E_2$ of zero and a reduction function $r:\tilde{U}_1\times \tilde{\Lambda}\mapsto \tilde{U}_2$ with the following properties. The reduction function $r$ is $k$ times continuously differentiable on $\tilde{U}_1\times \tilde{\Lambda}$ and the derivatives of $r$ are bounded and uniformly continuous on $\tilde{U}_1\times \tilde{\Lambda}$ with
\begin{equation*}
r(0,0)=0,\quad d_1r[0,0]=0.
\end{equation*}
The graph
\begin{equation*}
X_C^{\mu}=\{u_1+r(u_1,\mu)\in \tilde{U}_1\times \tilde{U}_2\ :\ u_1\in \tilde{U}_1\},
\end{equation*}
is a Hamiltonian center manifold for \eqref{fulleqthm} with the following properties:
\begin{itemize}
\item Through every point in $X_C^\mu$ there passes a unique solution of \eqref{fulleqthm} that remains on $X_C^\mu$ as long as it remains in $\tilde{U}_1\times \tilde{U}_2$. We say that $X_C^\mu$ is a locally invariant manifold of \eqref{fulleqthm}.
\item Every small bounded solution $u(x)$, $x\in \mathbb{R}$ of \eqref{fulleqthm} that satisfies $u_1(x)\in \tilde{U}_1$, $u_2(x)\in \tilde{U}_2$ lies completely in $X_C^\mu$.
\item Every solution $u_1$ of the reduced equation 
\begin{equation}
\dot{u}_1=K_1u_1+\mathcal{F}_1(u_1+r(u_1,\mu),\mu)\label{thmredeq},
\end{equation}
generates a solution
\begin{equation*}
u(x)=u_1(x)+r(u_1(x),\mu)
\end{equation*}
of \eqref{fulleqthm}.
\item $X_C^\mu$ is a symplectic submanifold of $E$, and the flow determined by the Hamiltonian system $(X_C^\mu,\Omega_C^\mu,H^\mu)$, where $\Omega_C^\mu$ is the reduced symplectic structure and $H^\mu$ the reduced Hamiltonian (see \eqref{defsym}--\eqref{defham}), coincides with the flow on $X_C^\mu$ determined by $(E,\Omega,\mathcal{H}^\mu)$. The reduced equation \eqref{thmredeq} represents Hamilton's equations for $(X_C^\mu,\Omega_C^\mu,H^\mu)$.
\item 
If \eqref{fulleqthm} is reversible, that is if there exists a linear symmetry $S$ which anticommutes with the right hand side of \eqref{fulleqthm}, then the reduction function $r$ can be chosen so that it commutes with $S$.
\end{itemize}
\end{theorem}
In our case we have $E=M$ and \eqref{fulleqthm} corresponds to \eqref{hameqlin2}.
We will use the same arguments as in \cite{BGT} when showing that hypothesis $H1-H4$ are satisfied.
Note that $H3$ is satisfied, by the following theorem:
\begin{lemma}
There exist constants $C,a_0>0$ such that 
\begin{equation}
\norm{(L-\i aI)^{-1}}_{M\mapsto M}\leq \frac{C}{\abs{a}}\label{resineq}
\end{equation}
for all $\abs{a}> a_0$.
\end{lemma}
The proof of this lemma is very similar to the proof of Lemma 3.4 in \cite{Groves2007} and will therefore be omitted. It follows from \eqref{linopk} that \eqref{resineq} holds for $K$ as well. In particular we get from \eqref{resineq} that the resolvent set of $K$ is nonempty, which implies that $K$ is closed.

Let $\i a$ be an element in the resolvent set of $K$. It follows from the Kondrachov embedding theorem that 
\begin{equation*}
K:\mathcal{D}(K)\mapsto M,
\end{equation*}
has compact resolvent. This implies that the spectrum of $K$ consists of an at most countable number of isolated eigenvalues with finite multiplicity. Combining this with the results of section \ref{spectrum} we can conclude that there exists $\xi>0$ such that
\begin{equation*}
\text{Spec}(K)=\{\lambda\in\text{Spec}(K)\ : \ \abs{\text{Re}(\lambda)}>\xi\}\cup \{\lambda \in \text{Spec}(K)\ :\ \text{Re}(\lambda)=0\},
\end{equation*}
that is, the part of the spectrum which lies on the imaginary axis is separated from the rest of the spectrum. This allows us to define the spectral projection $P$, corresponding to the imaginary part of the spectrum:
\begin{equation}
P=-\frac{1}{2\pi \i}\int_\gamma (K-\lambda)^{-1}d\lambda, \label{specprojection}
\end{equation}
where $\gamma$ is a curve surrounding the imaginary part of the spectrum and which lies in the resolvent set.

We check hypotheses $H1,H2$ and $H4$ of Theorem \ref{cmt}. From Lemma \ref{change_ofvariables_lemma} we get that $H4$ is satisfied, with $\Lambda:=\Delta$ and $U:=\tilde{V}\cap\mathcal{D}(K)$. 
Let $E_1=PM$ and let $E_2=(I-P)M$. It follows from Theorem 6.17 chapter III in \cite{KA}, together with the fact that the imaginary part of the spectrum of $K$ consists of a finite number of eigenvalues with finite multiplicity (see section \ref{spectrum}), that $H1$ and $H2$ are satisfied.

By the center-manifold theorem, there exist neighborhoods $\tilde{U}_1\subset U \cap E_1$, $\tilde{\Lambda}\subset\Lambda$ of zero and a reduction function $r:\tilde{U}_1\times \tilde{\Lambda}\mapsto\tilde{U}_2$ such that $r(0,0)=d_1r[0,0]=0$ and
\begin{equation*}
X_C^\mu=\{u_1+r(u_1,\mu)\ : \ u_1\in \tilde{U}_1\},
\end{equation*}
is a center manifold for \eqref{hameqlin2}. We then have the Hamiltonian system $(X_C^\mu,\Omega_C^\mu, H^\mu)$, where 
\begin{align}
(\Omega_C^\mu)_m(u_1,u_1^*)&=(\tilde{\Omega}^\mu)_{m+r(m,\mu)}(u_1+d_1r[m,\mu](u_1),u_1^*+d_1r[m,\mu](u_1^*)),\label{defsym}\\
H^\mu(u_1)&=\tilde{\mathcal{H}}^\mu(u_1+r(u_1,\mu)).\label{defham}
\end{align}
Thus far we have obtained a center manifold parametrized on the coordinate chart $\tilde{U}_1$ with coordinate map $\vartheta\ :\ X_C^\mu\mapsto \tilde{U}_1$ given by
\begin{equation*}
\vartheta^{-1}(u_1)=u_1+r(u_1,\mu).
\end{equation*}
Since the change of variables introduced in section \ref{changevar} is not explicit, we change parametrization by introducing the coordinate chart $\tilde{W}_1=(\mathrm{d}G^0[0])^{-1}(\tilde{U}_1)$ and the coordinate map $\varphi\ :\ X_C^\mu\mapsto \tilde{W}_1$ given by
\begin{equation*}
\varphi^{-1}(w_1)=w_1+h(w_1,\mu),
\end{equation*}
where $h\ :\ \tilde{W}_1\times \tilde{\Lambda}\mapsto V$ is a new reduction function defined by
\begin{equation*}
h(w_1,\mu)=(G^{\mu})^{-1}(\mathrm{d}G^0[0](w_1)+r(\mathrm{d}G^0[0](w_1),\mu))-w_1.
\end{equation*}
By construction $w_1+h(w_1,\mu)\in \mathcal{D}(v_{\mathcal{H}^\mu})$ which means that $X_C^\mu=\{w_1+h(w_1,\mu)\ : \ w_1\in \tilde{W}_1\}\subset \mathcal{D}(v_{\mathcal{H}^\mu})$, that is we have obtained a center manifold for \eqref{hamiltons_eq}. As before we have the Hamiltonian system $(X_C^\mu,\Omega_C^\mu,H^\mu)$, where
\begin{align*}
(\Omega_C^\mu)_m(w_1,w_1^*)&=\Omega(w_1+\mathrm{d}_1h[m,\mu](w_1),w_1^*+\mathrm{d}_1h[m,\mu](w_1^*)),\\
H^\mu(w_1)&=\mathcal{H}^\mu(w_1+h(w_1,\mu)).
\end{align*}
Note also that since $w_1+h(w_1,\mu)\in \mathcal{D}(v_{\mathcal{H}^\mu})$ it follows that 
\begin{equation*}
B_i^\mu(w_1+h(w_1,\mu))=0, i=0,1,\ y=0,1,
\end{equation*} 
for $w_1\in \tilde{W}_1$, where 
\begin{equation*}
B_1^\mu(u)=-\frac{\rho\phi_{1y}}{1-\eta}+\frac{\rho F_1(u,\mu)}{1-\eta},\ B_2^\mu(u)=-\frac{\phi_{2y}}{\eta+h}+\frac{F_2(u,\mu)}{\eta+h}.
\end{equation*}
Also note that 
\begin{equation}\label{boundary-prop}
\mathrm{d}\mathcal{H}^\mu[w](w^*)=\Omega(v_{\mathcal{H}^\mu}(w),w^*)+\int_0^{2\pi}B_1^\mu(w)\phi_1^{w^*}\vert_{y=1}\ \mathrm{d}z+\int_0^{2\pi}B^\mu_2(w)\phi_2^{w^*}\vert_{y=1}\ \mathrm{d}z,
\end{equation}
for $w\in N\cap \mathcal{D}(L)$, $w^*\in M$, which follows from the definition of $v_{\mathcal{H}^\mu}$ and $\mathcal{D}(v_{\mathcal{H}^\mu})$. 
\\
From Darboux's theorem (see \cite[Theorem 4]{BG}) there exists a near identity change of variables
\begin{equation*}
\hat{w}_1=w_1+\Upsilon(w_1,\mu),
\end{equation*}
such that $\Omega_C^\mu$ is transformed into $\Psi$, where 
\begin{equation*}
\Psi(w_1,w_1^*)=\Omega(w_1,w_1^*)
\end{equation*}
The coordinate map is then given by $\hat{w}_1+\hat{h}(\hat{w}_1,\mu)$, where $\hat{h}: \tilde{W}_1\times\tilde{\Lambda} \mapsto \tilde{W}_1\times \tilde{W}_2$ and $\hat{h}(0,0)=\mathrm{d}_1\hat{h}[0,0]=0$. In order to simplify the notation we immediately remove the $\wedge$ accent.\\ 
The next result is a generalization of \cite[Theorem 4.4]{Groves2017}.
\begin{theorem}\label{normalform}
Consider an $n+2$-degree of freedom Hamiltonian system
\begin{alignat}{3}
\dot{q}_i&=\frac{\partial{H}^\mu}{\partial p_i}, &&\dot{p}_i=-\frac{\partial H^\mu}{\partial q_i},\ i&&=1,2,\ldots ,n\label{abstractsys1}\\
\dot{q}_{n+1}&=\frac{\partial H^\mu}{\partial p_{n+1}},\  &&\dot{p}_{n+1}=-\frac{\partial H^\mu}{\partial q_{n+1}},&&\\
\dot{q}_{n+2}&=\frac{\partial H^\mu}{\partial p_{n+2}}, \ &&\dot{p}_{n+2}=-\frac{\partial H^\mu}{\partial q_{n+2}}&&\label{abstractsys3}
\end{alignat}
where $p_{n+1},\ p_{n+2}$ are cyclic variables so that $q_{n+1},\ q_{n+2}$ are conserved quantities and $\mu\in \mathbb{R}^l$ is a parameter. There exists a near-identity canonical change of variables $(q,p,q_{n+1},p_{n+1},q_{n+2},p_{n+2})\mapsto (Q,P,Q_{n+1},P_{n+1},Q_{n+2},P_{n+2})$ with the properties that $P_{n+1},\ P_{n+2}$ are cyclic, $Q_{n+1}=q_{n+1},\ Q_{n+2}=q_{n+2}$ and the lower order Hamiltonian system
\begin{equation*}
\dot{Q}_i=\frac{\partial H^\mu}{\partial P_i}(Q,P,0),\ \dot{P}_i=-\frac{\partial H^\mu}{\partial Q_i}(Q,P,0),\ i=1,2,\ldots, n,
\end{equation*}
adopts its usual normal form.  
\end{theorem} 
Theorem \ref{normalform} is useful when dealing with Hamiltonian systems with cyclic variables since it essentially allows us to treat the corresponding conserved quantities as parameters. We will apply Theorem \ref{normalform} in section \ref{Hamiltonian-Hopf bifurcation} when studying the Hamiltonian-Hopf bifurcation.
\section{Doubly periodic waves}\label{periodic solutions}
In this section we will examine solutions of \eqref{hamiltons_eq} that, in addition to being periodic in $z$, are periodic in $x$ with some period $P_x$. These correspond to doubly periodic solutions of the governing equations \eqref{goveqfinal1}--\eqref{goveqfinal8}.
Fix parameters $(\beta,\alpha,\theta_1,\theta_2)$ so that $\theta_1,\theta_2\neq \pm\pi/2$, $\cos^2(\theta_1)\leq\cos^2(\theta_2)$ and $(\beta,\alpha)\in I_{\rho,h,\theta_1}$. Then, according to the discussion in section \ref{mode1}, there exists $\nu_0$ such that the imaginary part of the spectrum of $L$ consists of $\{\pm \i\kappa_0,\pm \i s, 0\}$ where $\pm \i\kappa_0$ are algebraically simple mode $0$ eigenvalues with eigenvectors $v_{\kappa_0}^0$, $\bar{v}_{\kappa_0}^0$, $\pm \i s$ are algebraically simple mode $\pm 1$ eigenvalues with eigenvectors $\exp(\i z)v_s^1$, $\exp(-\i z)\bar{v}_s^1$ and $0$ is a geometrically double mode $\pm 1$ eigenvalue with eigenvectors $\exp(\i z)v_0^1$, $\exp(-\i z)\bar{v}_0^1$.
Moreover, $0$ is trivially a mode $0$ eigenvalue of algebraic multiplicity $4$, with eigenvectors $e_1,e_2$ and corresponding generalized eigenvectors $f_1, f_2$. In addition we assume that $\i s$ is nonresonant with $\i\kappa_0$. In addition of considering bifurcations around $v_0$ we will also consider bifurcations in $2\pi/P_x$. In anticipation of this we let $\mu=\mu_1$ throughout this section.
Let
\begin{align*}
&V_{\kappa_0}^0=\frac{v_{\kappa_0}^0}{\sqrt{\abs{c_1}}},\ V_0^1=\frac{v_0^1}{\sqrt{\abs{c_2}}},\ V_s^1=\frac{v_s^1}{\sqrt{\abs{c_3}}},\\
&\tilde{e}_1=\frac{-\frac{1}{\rho\cos^2(\theta_1)}(h\alpha-\cos^2(\theta_1))e_1+e_2}{\sqrt{c_4}},\ \tilde{e}_2=\frac{e_1-\frac{1}{\cos^2(\theta_1)}(\alpha-\rho\cos^2(\theta_1))e_2}{\sqrt{c_4}},\\
& \tilde{f}_i=\frac{f_i}{\sqrt{c_4}},\ i=1,2,
\end{align*} 
where
\begin{align*}
c_1&=\frac{4\pi\kappa_0}{\cos^2(\theta_1)}(\beta_0^*(\kappa_0)-\beta),\\
c_2&=
\begin{cases}
&\frac{4\pi\cos(\theta_1-\theta_2)\nu_0}{\cos^2(\theta_2)}(\beta_1^*(0)-\beta),\text{ if }\cos(\theta_1-\theta_2)\neq 0\\
&\frac{4\pi\nu_0\cos(\theta_2)}{\cos(\theta_2)}\left(\frac{\rho}{\tanh(\nu_0)}+\frac{1}{\tanh(h\nu_0)}\right), \text{ if }\cos(\theta_1-\theta_2)=0,
\end{cases}\\
c_3&=\begin{cases}
&\frac{4\pi\tilde{\gamma}_{1}^2(\nu_0\cos(\theta_1-\theta_2)+s)}{(\nu_0\cos(\theta_2)+s\cos(\theta_1))^2}(\beta_1^*(s)-\beta),\text{ if } s+\nu_0\cos(\theta_1-\theta_2)\neq 0\\
&\frac{4\pi \tilde{\gamma}_{1}\cos(\theta_1)}{\nu_0\cos(\theta_2)+s\cos(\theta_1)}\left(\frac{\rho}{\tanh(\tilde{\gamma}_{1})}+\frac{1}{\tanh(h\tilde{\gamma}_{1})}\right),\text{ if } s+\nu_0\cos(\theta_1-\theta_2)= 0,
\end{cases}\\
c_4&=\frac{2\pi h}{\cos^2(\theta_1)}\left(\cos^2(\theta_1)\left(\rho+\frac{1}{h}\right)-\alpha\right).
\end{align*}
Then 
\begin{align*}
&\Omega(V_{\kappa_0}^0,\bar{V}_{\kappa_0}^0)=\text{sgn}(c_1)\i,&& \Omega(\exp(\i z)V_0^1,\exp(-\i z)\bar{V}_0^1)=\text{sgn}(c_2)\i,\\
&\Omega(\exp(\i z)V_s^1,\exp(-\i z)\bar{V}_s^1)=\text{sgn}(c_3)\i,&& \Omega(\tilde{e}_i,\tilde{f}_i)=1,\ i=1,2,
\end{align*}
and all other combinations are equal to $0$. The signs of the coefficients will not affect the subsequent analysis, so we assume for definiteness that
\begin{equation*}
\text{sgn}(c_1)=\text{sgn}(c_2)=\text{sgn}(c_3)=1.
\end{equation*}
This can for example be achieved by choosing $\beta<\min(\beta_0^*(\kappa_0),\beta_1^*(0),\beta_1^*(s))$ and $\theta_1-\theta_2\in (-\pi/2,\pi/2)$. Hence, $\{V_{\kappa_0}^0,\bar{V}_{\kappa_0}^0,\exp(\i z)V_0^1,\exp(-\i z)\bar{V}_0^1,\exp(\i z)V_s^1,\exp(-\i z)\bar{V}_s^1,\tilde{e}_1,\tilde{e}_2,\tilde{f}_1,\tilde{f}_2\}$ is a symplectic basis of $(E_1,\Psi)$. We introduce coordinates on $X_C^{\mu_1}$ by writing
\begin{align*}
w_1&=(A,B,C,\bar{A},\bar{B},\bar{C},q_1,p_1,q_2,p_2)\\
&=AV_{\kappa_0}^0+B\exp(\i z)V_0^1+C\exp(\i z)V_s^1+\bar{A}\bar{V}_{\kappa_0}^0+\bar{B}\exp(-\i z)\bar{V}_0^1+\bar{C}\exp(-\i z)\bar{V}_s^1\\
&\quad+q_1\tilde{e}_1+p_1\tilde{f}_1+q_2\tilde{e}_2+p_2\tilde{f}_2.
\end{align*}
By construction 
\begin{equation*}
\Psi(w_1,w_1^*)=\i (A\bar{A}^*+B\bar{B}^*+C\bar{C}^*-\bar{A}A^*-\bar{B}B^*-\bar{C}C^*)+q_1p_1^*+q_2p_2^*-q_1^*p_1-q_2^*p_2.
\end{equation*}
Recall that $\mathcal{H}^{\mu_1}$ is invariant under the transformations $\phi_i\mapsto \phi_i+b_i$, $i=1,2$, for arbitrary $b_i\in \mathbb{R}$. This symmetry is inherited by the reduced system: the variables $q_1$ and $q_2$ are cyclic, that is $H^{\mu_1}$ is independent of $q_1,q_2$, which implies that $p_1,p_2$ are conserved. We may therefore set them to $0$ and recover the variables $q_1,q_2$ by quadrature. Due to this we introduce $\tilde{w}_1=(A,B,C,\bar{A},\bar{B},\bar{C})\in \widehat{E}_1\cap \tilde{W}_1$, where $\widehat{E}_1=\text{span}_\mathbb{C}\{V_{\kappa_0}^0,\exp(\i z)V_0^1,\exp(\i z)V_s^1,\bar{V}_{\kappa_0}^0,\exp(-\i z)\bar{V}_0^1,\exp(-\i z)\bar{V}_s^1\}$, and write Hamilton's equations for the reduced Hamiltonian system $(X_C^{\mu_1},\Psi,H^{\mu_1})$ as
\begin{equation}\label{redhameq1}
\dot{\tilde{w}}_1-v_{H^{\mu_1}}(\tilde{w}_1)=0,\ \tilde{\omega}_1\in \widehat{E}_1\cap\tilde{W}_1,
\end{equation}
where 
\begin{equation*}
v_{H^{\mu_1}}(\tilde{w}_1)=-\i\left(\begin{array}{c}
\frac{\partial H^{\mu_1}}{\partial\bar{A}}\\
\frac{\partial H^{\mu_1}}{\partial\bar{B}}\\
\frac{\partial H^{\mu_1}}{\partial\bar{C}}\\
-\frac{\partial H^{\mu_1}}{\partial A}\\
-\frac{\partial H^{\mu_1}}{\partial B}\\
-\frac{\partial H^{\mu_1}}{\partial C}
\end{array}
\right).
\end{equation*}
We next define a symplectic structure on $\widehat{E}_1$, by
\begin{equation*}
J(\tilde{w}_1)=J(A,B,C,\bar{A},\bar{B},\bar{C})=-\i (A,B,C,-\bar{A},-\bar{B},-\bar{C}),
\end{equation*}
and an inner product
\begin{equation*}
\langle \tilde{w}_1,\tilde{w}_1^*\rangle=A\bar{A}^*+B\bar{B}^*+C\bar{C}^*+\bar{A}A^*+\bar{B}B^*+\bar{C}C^*,
\end{equation*}
so that $v_{H^{\mu_1}}(\tilde{w}_1)=J\nabla H^{\mu_1}(\tilde{w}_1)$, and $\Psi(\tilde{w}_1,\tilde{w}_1^*)=-\langle J\tilde{w}_1,\tilde{w}_1^*\rangle$.
The next step is to normalize the period in $x$, which introduces the parameter $\kappa=2\pi/P_x$. We will consider values of $\kappa$ close to $\kappa_0$ and so we introduce a bifurcation parameter $\mu_2$ by writing $\kappa=\kappa_0+\mu_2$. Equation \eqref{redhameq1} can then be written as
\begin{equation}\label{redhameq2}
(\kappa_0+\mu_2)J\dot{\tilde{w}}_1+\nabla H^{\mu_1}(\tilde{w}_1)=0,\ \tilde{w}_1\in \widehat{E}_1\cap\tilde{W}_1, \ (\mu_1,\mu_2)\in \widehat{\Lambda},
\end{equation}
where $\widehat{\Lambda}$ is some neighborhood of the origin in $\mathbb{R}^2$. 
Next we want to consider solutions of \eqref{redhameq1} as elements of the Sobolev space $\mathcal{Y}:=H_{per}^1((0,2\pi),\widehat{E}_1)$, equipped with the norm $\norm{\tilde{w}_1}_{\mathcal{Y}}=(\norm{\tilde{w}_1}_{L_{per}^2((0,2\pi),\widehat{E}_1)}^2+\norm{\dot{\tilde{w}}_1}_{L_{per}^2((0,2\pi),\widehat{E}_1)}^2)^{1/2}$, where $\norm{\cdot}_{L_{per}^2((0,2\pi),\widehat{E}_1)}$ is the norm coming from the inner product
\begin{equation*}
(\tilde{w}_1,\tilde{w}_1^*)=\int_0^{2\pi}\langle\tilde{w}_1,\tilde{w}_1^*\rangle \ \mathrm{d}x.
\end{equation*}
Define
\begin{equation*}
T(\tilde{w}_1,\mu_1,\mu_2)=(\kappa_0+\mu_2)J\dot{\tilde{w}}_1+\nabla H^{\mu_1}(\tilde{w}_1),\ \tilde{w}_1\in B_R(0), \ (\mu_1,\mu_2)\in \widehat{\Lambda},
\end{equation*}
where $B_R(0)$ is the ball of radius $R$ centered at the origin in $\mathcal{Y}$ and where $R$ is chosen small enough so that  $\tilde{w}_1(x)\in \tilde{W}_1$, for $\tilde{w}_1\in B_R(0)$. Then \eqref{redhameq2} can be written
\begin{equation}\label{redhameq3}
T(\tilde{w}_1,\mu_1,\mu_2)=0.
\end{equation}
Equation \eqref{redhameq3} can be seen as the Euler-Lagrange equation of the action integral
\begin{equation*}
\mathcal{E}(\tilde{w}_1,\mu_1,\mu_2)=\int_0^{2\pi} \frac{\kappa_0+\mu_2}{2}\langle J\dot{\tilde{w}}_1,\tilde{w}_1\rangle+\tilde{H}^{\mu_1}(\tilde{w}_1)\ \mathrm{d}x,
\end{equation*}
with respect to the inner product $(\cdot,\cdot)$, that is
\begin{equation}\label{def-EL}
\mathrm{d}\mathcal{E}[\tilde{w}_1,\mu_1,\mu_2](\tilde{w}_1^*)=(T(\tilde{w}_1,\mu_1,\mu_2),\tilde{w}_1^*),\ \tilde{w}_1\in B_{R}(0),\ \tilde{w}_1^*\in \mathcal{Y},\ (\mu_1,\mu_2)\in\widehat{\Lambda}.
\end{equation}
We will find critical points of $\mathcal{E}$ by using a variational Lyapunov-Schmidt reduction; see for example \cite{Berti2007}. First note that
\begin{equation*}
\mathcal{K}:=\text{ker}(\mathrm{d}T[0,0,0])=\text{span}_\mathbb{C}\{\exp(\i x)V_{\kappa_0}^0,\exp(\i z)V_0^1,\exp(-\i x)\bar{V}_{\kappa_0}^0,\exp(-\i z)\bar{V}_0^1\}.
\end{equation*}
Next decompose
\begin{equation*}
\mathcal{Y}=\mathcal{K}\oplus\mathcal{K}^\perp,
\end{equation*}
and write $\tilde{w}_1=\zeta+\xi$, where $\zeta\in\mathcal{K}$, $\xi\in\mathcal{K}^\perp$. Let $\Pi$ be the projection onto $\mathcal{K}$, so that equation \eqref{redhameq3} can be decomposed as
\begin{align}
\Pi(T(\zeta+\xi,\mu_1,\mu_2))&=0,\label{bifeq}\\
(I-\Pi)(T(\zeta+\xi,\mu_1,\mu_2))&=0.\label{rangeeq}
\end{align}
Equation \eqref{rangeeq} can be solved using the implicit function theorem, which yields solutions of the form $\zeta+h_{LS}(\zeta,\mu_1,\mu_2)$, where $h_{LS}\ :\ \mathcal{K}_0\times \widehat{\Lambda}_0 \mapsto \mathcal{K}^\perp$ and $\mathcal{K}_0$, $\widehat{\Lambda}_0$ are open neighborhoods of the origin in $\mathcal{K}$ and $\mathbb{R}^2$ respectively. In particular we can assume that $\zeta+h_{LS}(\zeta,\mu_1,\mu_2)\subseteq B_R(0)$, for all $\zeta\in\mathcal{K}_0,\ (\mu_1,\mu_2)\in\widehat{\Lambda}_0$. In order to solve equation \eqref{bifeq} we define the reduced functional 
\begin{equation*}
\mathcal{E}_{LS}(\zeta,\mu_1,\mu_2)=\mathcal{E}(\zeta+h_{LS}(\zeta,\mu_1,\mu_2),\mu_1,\mu_2).
\end{equation*}
Note that for all $\zeta\in\mathcal{K}_0$ and $\zeta^*\in \mathcal{K}$,
\begin{align}
\mathrm{d}\mathcal{E}_{LS}[\zeta,\mu_1,\mu_2](\zeta^*)&=\mathrm{d}\mathcal{E}[\zeta+h_{LS}(\zeta,\mu_1,\mu_2),\mu_1,\mu_2](\zeta^*+\mathrm{d}h_{LS}[\zeta,\mu_1,\mu_2](\zeta^*))\nonumber\\
&=( T(\zeta+h_{LS}(\zeta,\mu_1,\mu_2),\mu_1,\mu_2),\zeta^*+\mathrm{d}h_{LS}[\zeta,\mu_1,\mu_2](\zeta^*))\nonumber\\
&=(\Pi T(\zeta+h_{LS}(\zeta,\mu_1,\mu_2),\mu_1,\mu_2),\zeta^*)\label{redvarlyap},
\end{align}
where we used that $\mathrm{d}_1h_{LS}[\zeta,\mu_1,\mu_2](\zeta^*)\in \mathcal{K}^\perp$ and that \eqref{rangeeq} is satisfied. The calculation \eqref{redvarlyap} shows that \eqref{bifeq} is the Euler-Lagrange equation of the action integral $\mathcal{E}_{LS}$, so solutions of \eqref{bifeq} are critical points of $\mathcal{E}_{LS}$. In order to find critical points of the functional $\mathcal{E}_{LS}$ we introduce coordinates in $\mathcal{K}$:
\begin{equation*}
\zeta=(A,B,\bar{A},\bar{B})=A\exp(\i x)V_{\kappa_0}^0+B\exp(\i z)V_0^1+\bar{A}\exp(-\i x)\bar{V}_{\kappa_0}^0+\bar{B}\exp(-\i z)\bar{V}_0^1,
\end{equation*}
and write $\mathcal{E}_{LS}(\zeta,\mu_1,\mu_2)=\mathcal{E}_{LS}(A,B,\bar{A},\bar{B},\mu_1,\mu_2)$. Then $\zeta$ is a critical point of $\mathcal{E}_{LS}$ if and only if $\nabla\mathcal{E}_{LS}(A,B,\bar{A},\bar{B},\mu_1,\mu_1)=0$, which is equivalent with
\begin{equation}\label{lyap1}
\begin{dcases}
\frac{\partial \mathcal{E}_{LS}}{\partial \bar{A}}(A,B,\bar{A},\bar{B},\mu_1,\mu_2)&=0,\\
\frac{\partial \mathcal{E}_{LS}}{\partial \bar{B}}(A,B,\bar{A},\bar{B},\mu_1,\mu_2)&=0.
\end{dcases}
\end{equation}
Recall that $\mathcal{H}^{\mu_1}$ is invariant under the transformation $z\mapsto z+z_0$, $z_0\in[0,2\pi]$. The reduction function $h$ can be chosen in such a way that it commutes with this transformation, which implies that the reduced Hamiltonian $H^{\mu_1}$ is invariant under the same transformation. Clearly the same is then true for $\mathcal{E}$, and in addition $\mathcal{E}$ is invariant under $x\mapsto x+x_0$, $x_0\in[0,2\pi]$. It follows that the reduced functional $\mathcal{E}_{LS}$ is invariant under rotations in both $x$ and $z$ as well. In terms of coordinates, this means that $\mathcal{E}_{LS}$ is invariant under the transformations
\begin{align*}
(A,B,\bar{A},\bar{B})&\mapsto (A,\exp(\i z_0)B,\bar{A},\exp(-\i z_0)\bar{B}),\\
(A,B,\bar{A},\bar{B})&\mapsto (\exp(\i x_0)A,B,\exp(-\i x_0)\bar{A},\bar{B}),
\end{align*}
 which implies that (see \cite[Sect VI, Lemma 2.1]{gols85})
\begin{equation*}
\mathcal{E}_{LS}(A,B,\bar{A},\bar{B},\mu_1,\mu_2)=\mathcal{E}_{LS}(\abs{A}^2,\abs{B}^2,\mu_1,\mu_2),
\end{equation*}
and so there exist functions $\Theta_i(\abs{A}^2,\abs{B}^2,\mu_1,\mu_2)$, $i=1,2$, such that
\begin{align*}
\frac{\partial \mathcal{E}_{LS}}{\partial \bar{A}}(A,B,\bar{A},\bar{B},\mu_1,\mu_2)&=A\Theta_1(\abs{A}^2,\abs{B}^2,\mu_1,\mu_2),\\
\frac{\partial \mathcal{E}_{LS}}{\partial \bar{B}}(A,B,\bar{A},\bar{B},\mu_1,\mu_2)&=B\Theta_2(\abs{A}^2,\abs{B}^2,\mu_1,\mu_2).
\end{align*}
The system \eqref{lyap1} becomes
\begin{equation}\label{lyap2} 
\begin{dcases}
&A\Theta_1(\abs{A}^2,\abs{B}^2,\mu_1,\mu_2)=0,\\
&B\Theta_2(\abs{A}^2,\abs{B}^2,\mu_1,\mu_2)=0.
\end{dcases}
\end{equation}
It is clear that $A=B=0$ is a solution of \eqref{lyap2}, for all $\mu_1,\mu_2$. In order to find nontrivial solutions we apply the implicit function theorem, and therefore want to show that 
\begin{equation*}
\text{det}\left(\begin{array}{cc}
\frac{\partial \Theta_1}{\partial\mu_1}(0,0,0,0)& \frac{\partial \Theta_1}{\partial\mu_2}(0,0,0,0)\\
\frac{\partial \Theta_2}{\partial\mu_1}(0,0,0,0)& \frac{\partial \Theta_2}{\partial\mu_2}(0,0,0,0)
\end{array}\right)\neq 0.
\end{equation*}
Denote by $\mathcal{E}_{LS,i}^{mn}$ the part of $\mathcal{E}_{LS}$ that is homogeneous of order $i$ in $(\abs{A}^2,\abs{B}^2)$, $m$ in $\mu_1$ and $n$ in $\mu_2$. In the same way we denote by $h_{LS,ijkl}^{mn}$, $h_{ij0kl0}^{m}$ the parts of $h_{LS}$, $h$ which are homogeneous of order $i,j,k,l,m,n$ in $A,B,\bar{A},\bar{B},\mu_1,\mu_2$, respectively. Then
\begin{align*}
\mathcal{E}_{LS,1}^{10}(\abs{A}^2,\abs{B}^2,\mu_1,\mu_2)&=d_1^{10}\mu_1\abs{A}^2+d_2^{10}\mu_1\abs{B}^2,\\
\mathcal{E}_{LS,1}^{01}(\abs{A}^2,\abs{B}^2,\mu_1,\mu_2)&=d_1^{01}\mu_2\abs{A}^2+d_2^{01}\mu_2\abs{B}^2,
\end{align*}
and
\begin{align*}
&\frac{\partial \Theta_1}{\partial\mu_1}(0,0,0,0)=d_1^{10},\quad \frac{\partial \Theta_1}{\partial\mu_2}(0,0,0,0)=d_1^{01},\\
&\frac{\partial \Theta_2}{\partial\mu_1}(0,0,0,0)=d_2^{10},\quad \frac{\partial \Theta_2}{\partial\mu_2}(0,0,0,0)=d_2^{01}.
\end{align*}
First note that 
\begin{align}
d_1^{01}&=\mathrm{d}^2\mathrm{d}_{\mu_2}\mathcal{E}[0,0,0](\exp(\i x)V_{\kappa_0}^0,\exp(-\i x)\bar{V}_{\kappa_0}^0)+\mathrm{d}^2\mathcal{E}[0,0,0](\exp(\i x)V_{\kappa_0}^0,h_{LS,0010}^{01})\nonumber\\
&\quad+\mathrm{d}^2\mathcal{E}[0,0,0](\exp(-\i x)\bar{V}_{\kappa_0}^0,h_{LS,1000}^{01})\label{d1formula},
\end{align}
and since \eqref{redhameq3} is the Euler-Lagrange equation of $\mathcal{E}$, we have that 
\begin{equation*}
\mathrm{d}^2\mathcal{E}[0,0,0](\exp(\i x)V_{\kappa_0}^0,h_{LS,0010}^{01})=(\mathrm{d}T[0,0,0](\exp(\i x)V_{\kappa_0}^0),h_{LS,0010}^{01})=0,
\end{equation*}
since $\exp(\i x)V_{\kappa_0}^0\in\mathcal{K}$. 
In the same way we have that $\mathrm{d}^2\mathcal{E}[0,0,0](\exp(-\i x)\bar{V}_{\kappa_0}^0,h_{LS,1000}^{01})=0$. Moreover, from the definition of $\mathcal{E}$ we find that
\begin{align*}
\mathrm{d}^2\mathrm{d}_{\mu_2}\mathcal{E}[0,0,0](\exp(\i x)V_{\kappa_0}^0,\exp(-\i x)\bar{V}_{\kappa_0}^0)&=\i \int_0^{2\pi}(J(\exp(\i x)V_{\kappa_0}^0),\exp(-\i x)\bar{V}_{\kappa_0}^0)\ dx\\
&=\i \int_0^{2\pi}-\Psi(\exp(\i x)V_{\kappa_0}^0,\exp(-\i x)\bar{V}_{\kappa_0}^0)\ dx\\
&=2\pi.
\end{align*}
Hence, we get from \eqref{d1formula} that $d_1^{01}=2\pi$. A similar calculation shows that $d_2^{01}=0$.
Next, using the same methods as above we find that
\begin{align}
d_2^{10}&=\mathrm{d}^2\mathrm{d}_{\mu_1}\mathcal{E}[0,0,0](\exp(\i z)V_0^1,\exp(-\i z)\bar{V}_0^1)\nonumber\\
&=\mathrm{d}^2\mathrm{d}_{\mu_1}\mathcal{H}^0[0](\exp(\i z)V_0^1,\exp(-\i z)\bar{V}_0^1)+d^2\mathcal{H}^0[0](\exp(\i z)V_0^1,h_{000010}^{1})\nonumber\\
&\quad+\mathrm{d}^2\mathcal{H}^0[0](\exp(-iz)\bar{V}_0^1,h_{010000}^{1}),\label{c2formula}
\end{align}
and
\begin{align*}
\mathrm{d}^2\mathcal{H}^0[0](\exp(\i z)V_0^1,h_{000010}^{1})&=\Omega(L\exp(\i z)V_0^1,h_{000010}^{1})+\int_0^{2\pi}\mathrm{d}B_1^0[0](\exp(\i z)V_0^1)\phi_1^{h_{000010}^{1}}\vert_{y=1}\ \mathrm{d}z\\
&\quad+ \int_0^{2\pi}\mathrm{d}B_2^0[0](\exp(\i z)V_0^1)\phi_2^{h_{000010}^{1}}\vert_{y=1}\ \mathrm{d}z\\
&=0,
\end{align*}
where we used that $\exp(\i z)V_0^1\in \mathcal{D}(L)$ and $L\exp(\i z)V_0^1=0$. In the same way we find that $\mathrm{d}^2\mathcal{H}^0[0](\exp(-\i z)\bar{V}_0^1,h_{010000}^{1})=0$. Equation \eqref{c2formula} then tells us that
\begin{align*}
d_2^{10}&=\mathrm{d}^2\mathrm{d}_{\mu_1}\mathcal{H}^0[0](\exp(\i z)V_0^1,\exp(-\i z)\bar{V}_0^1)\\
&=\frac{4\pi\nu_0}{c_2\cos^2(\theta_2)}\Bigg[\frac{\cos^2(\theta_2)}{2\nu_0}\bigg(\frac{\rho}{\tanh(\nu_0)}+\frac{1}{\tanh(h\nu_0)}\bigg)-\frac{\cos^2(\theta_2)}{2}\bigg(\frac{\rho}{\sinh^2(\nu_0)}+\frac{h}{\sinh^2(h\nu_0)}\bigg)-\beta\Bigg]
\end{align*}
and this is nonzero precisely when $\beta\neq  \tilde{\beta}(\nu_0)$. This condition is automatically fulfilled in our case, since we assume that $(\beta,\alpha)\in I_{\rho,h,\theta_1}$ with $\cos^2(\theta_1)\leq \cos^2(\theta_2)$. Moreover
\begin{equation*}
\text{det}\left(\begin{array}{cc}
\frac{\partial \Theta_1}{\partial\mu_1}(0,0,0,0)& \frac{\partial \Theta_1}{\partial\mu_2}(0,0,0,0)\\
\frac{\partial \Theta_2}{\partial\mu_1}(0,0,0,0)& \frac{\partial \Theta_2}{\partial\mu_2}(0,0,0,0)
\end{array}\right)=-2\pi d_2^{10}\neq 0.
\end{equation*}
It now follows from the implicit function theorem that there exist nontrivial solutions of \eqref{lyap2}, and in conclusion we have the following result.
\begin{theorem}\label{doublyperiodic_thm}
If $\cos^2(\theta_1)\leq\cos^2(\theta_2)$, $\theta_1,\theta_2\neq\pm\pi/2$ and $(\beta,\alpha)\in I_{\rho,h,\theta_1}$, there exist $\kappa_0,\nu_0>0$, $\epsilon_1,\epsilon_2>0$ and functions $\mu_i\ : \ (0,\epsilon_i)\mapsto (0,\infty)$, $i=1,2$ such that
\begin{equation*}
w_1=\tilde{w}_1(A,B,\bar{A},\bar{B})+h(\tilde{w}_1(A,B,\bar{A},\bar{B}),\mu_1(\abs{A}^2,\abs{B}^2)),
\end{equation*}
is a doubly periodic travelling wave, with periods $2\pi/(\kappa_0+\mu_2(\abs{A}^2,\abs{B}^2))$ in $x$, $2\pi/(\nu_0+\mu_1(\abs{A}^2,\abs{B}^2))$ in $z$, for all $A,B$ such that $\abs{A}^2<\epsilon_1$, $\abs{B}^2<\epsilon_2$, where
\begin{equation}\label{specialsolutions}
\tilde{w}_1(A,B,\bar{A},\bar{B})=\check{w}_1+h_{LS}(\check{w}_1,\mu_1(\abs{A}^2,\abs{B}^2),\mu_2(\abs{A}^2,\abs{B}^2)),
\end{equation}
with
\begin{equation*}
\check{w}_1=A\exp(\i (\kappa_0+\mu_2)x)V_{\kappa_0}^0+B\exp(\i (\nu_0+\mu_1)z)V_0^1+\bar{A}\exp(-\i (\kappa_0+\mu_2)x)\bar{V}_{\kappa_0}^0+\bar{B}\exp(-\i(\nu_0+\mu_1) z)\bar{V}_0^1.
\end{equation*}
If instead $\cos^2(\theta_2)\leq\cos^2(\theta_1)$, $\theta_1,\theta_2\neq\pm\pi/2$, the theorem still holds with $(\beta,\alpha)\in I_{\rho,h,\theta_2}$.
\end{theorem}
\begin{remark}
As indicated in section \ref{mode1} we could allow for $L$ to have additional mode $k$ eigenvalues, as long as they are nonresonant with $\kappa_0$. In this case the corresponding eigenvectors will not be in the kernel of $\mathrm{d}T[0,0,0]$ and will therefore not affect the calculations once the Lyapunov-Schmidt reduction is carried out. It is also possible to obtain a similar result when $\cos^2(\theta_1)\leq \cos^2(\theta_2)$, $(\beta,\alpha)\in II_{\rho,h,\theta_1}$,  or $\cos^2(\theta_2)\leq \cos^2(\theta_1)$, $(\beta,\alpha)\in I_{\rho,h,\theta_2}\cup II_{\rho,h,\theta_2}$. In these cases there could possibly be some additional mode $0$ eigenvalues $\pm \i\kappa_1$. However, as explained above, as long as $\i \kappa_1$ is nonresonat with $\i \kappa_0$, they will have no impact on the calculations.
\end{remark}
\begin{remark}
Since we in particular need to assume in Theorem \ref{doublyperiodic_thm} that $\theta_2\neq \pm \pi/2$, we cannot directly obtain waves that are periodic in $Z$ with a bounded profile in the direction $X$. However, such solutions can be obtained using the Lyapunov-center theorem as in \cite[Theorem 3.9]{Groves2007}. Similarly, to obtain waves that are periodic in the direction $X$ with a bounded profile in $Z$, we could again apply the Lyapunov-center theorem as in \cite[Theorem 5]{Groves2001}.
\end{remark}

\section{Hamiltonian-Hopf bifurcation}\label{Hamiltonian-Hopf bifurcation}
\begin{figure}[H]
\centering
\includegraphics{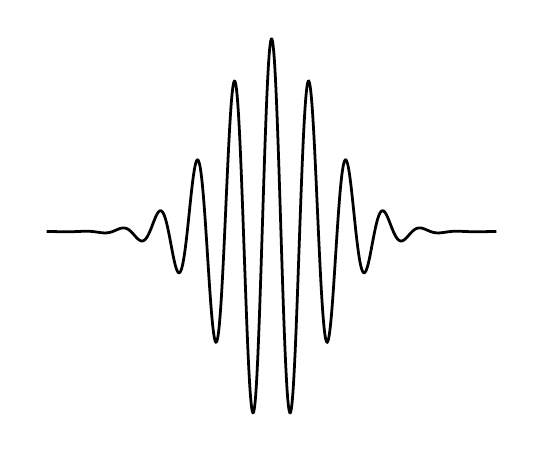}
\includegraphics{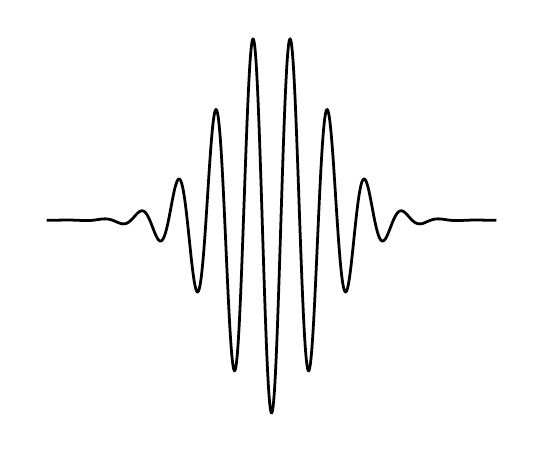}
\includegraphics{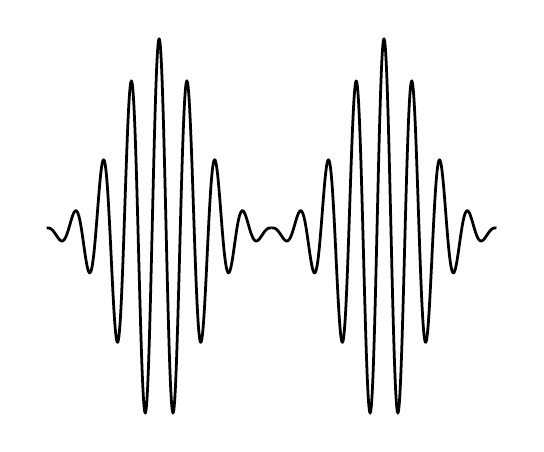}
\includegraphics{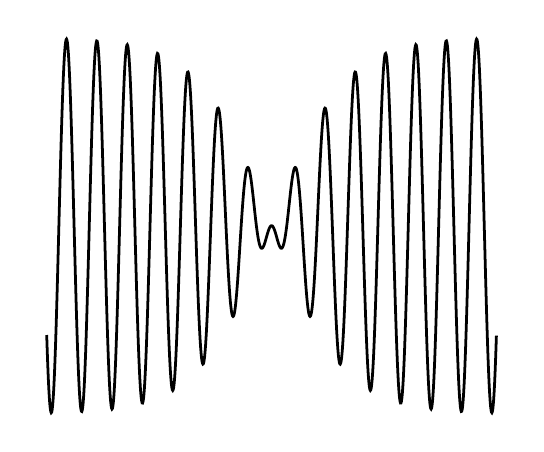}
\caption{Typical wave profiles of the different solutions found from theorem \ref{darkbright}. From left to right: Bright solitary wave of elevation, bright solitary wave of depression, multipulse bright solitary wave of elevation, dark solitary wave.}
\label{darkbrightwaveprofiles}
\end{figure}

In this section we consider the Hamiltonian-Hopf bifurcation occurring at some critical value $\nu_0$, which was discussed in section \ref{mode1}. For definiteness we will focus on the case when $\theta_1\neq \pm\pi/2$, so that $\i s$ is a mode $k$ eigenvalue of algebraic multiplicity $2$ if and only if $(\beta,\alpha)=(\beta_1^*(s),\alpha_1^*(s))$. In addition we must then have that $s+\nu_0\cos(\theta_1-\theta_2)\neq 0$. We therefore fix parameters $(\beta,\alpha,\theta_1,\theta_2,\nu_0)$, where $(\beta,\alpha)=(\beta_1^*(s),\alpha_1^*(s))$, $\theta_1\neq 0,\ \theta_1\neq \pm\pi/2$ and $(\beta,\alpha)\in III_{\rho,h,\theta_1}\cup IV_{\rho,h,\theta_1}$. Then the spectrum of $L$ consists of $\{\pm \i s,0\}$, where $\pm \i s$ are mode $\pm 1$ eigenvalues of algebraic multiplicity $2$ with eigenvectors $\exp(\i z)v_s^1,\ \exp(-\i z)\bar{v}_s^1$ and corresponding generalized eigenvectors $\exp(\i z)u_s^1$, $\exp(-\i z)\bar{u}_s^1$. Recall again that $0$ is an eigenvalue of algebraic multiplicity $4$ with eigenvectors $e_1,\ e_2$ and corresponding generalized eigenvectors $f_1$, $f_2$. Let
\begin{align*}
V_s^1&=\frac{v_s^1}{\sqrt{\tau_1}},\ U_s^1=\frac{1}{\sqrt{\tau_1}}\bigg(u_s^1-\frac{\i\tau_2v_s^1}{2\tau_1}\bigg),\\ 
\tilde{e}_1&=\frac{\frac{1}{\rho\cos^2(\theta_1)}(h\alpha-\cos^2(\theta_1))e_1-e_2}{\sqrt{\tau_3}},\ \tilde{e}_2=\frac{-e_1+\frac{1}{\cos^2(\theta_1)}(\alpha-\rho\cos^2(\theta_1))e_2}{\sqrt{\tau_3}},\\
\tilde{f}_i&=\frac{f_i}{\sqrt{\tau_3}},\ i=1,2,
\end{align*} 
where 
\begin{align*}
\tau_1&=\Omega(\exp(\i z)v_s^1,\exp(-\i z)\bar{u}_s^1),\ \i\tau_2=\Omega(\exp(\i z)u_s^1,\exp(-\i z)\bar{u}_s^1),\\
\tau_3&=\frac{2\pi h}{\cos^2(\theta_1)}\left(\alpha-\cos^2(\theta_1)\left(\rho+\frac{1}{h}\right)\right).
\end{align*}
We find in particular that
\begin{equation*}
\Omega(\exp(\i z)v_s^1,\exp(-\i z)\bar{u}_s^1)=-\frac{2\pi(s+\nu_0\cos(\theta_1-\theta_2))\tilde{\gamma}_1^2}{(s\cos(\theta_1)+\nu_0\cos(\theta_2))^2}\frac{\mathrm{d\beta_1^*(s)}}{\mathrm{d}s}.
\end{equation*}
Since we assume that $\i s$ is of algebraic multiplicity $2$ we have in particular that $\tau_1\neq 0$, and we assume for definiteness that $\tau_1>0$. This is for example achieved when $\frac{\mathrm{d\beta_1^*(s)}}{\mathrm{d}s}<0$, $s+\nu_0\cos(\theta_1-\theta_2)>0$. Then
\begin{equation*}
\Omega(\exp(\i z)V_s^1,\exp(-\i z)\bar{U}_s^1)=\Omega(\exp(-\i z)\bar{V}_s^1,\exp(\i z)U_s^1)=1, \Omega(\tilde{e}_i,\tilde{f}_i)=1,\ i=1,2,
\end{equation*}
and all other combinations are equal to zero. Hence, the set of vectors
\begin{equation*}
\{\exp(\i z)V_s^1,\exp(\i z)U_s^1,\exp(-\i z)\bar{V}_s^1,\exp(-\i z)\bar{U}_s^1,\tilde{e}_1,\tilde{e}_2,\tilde{f}_1,\tilde{f}_2\},
\end{equation*}
is a symplectic basis of $(E_1,\Psi)$. We introduce coordinates on $X_C^\mu$ by writing
\begin{equation*}
w_1=A\exp(\i z)V_s^1+B\exp(\i z)U_s^1+\bar{A}\exp(-\i z)\bar{V}_s^1+\bar{B}\exp(-\i z)\bar{U}_s^1+q_1\tilde{e}_1+p_1\tilde{f}_1+q_2\tilde{e}_2+p_2\tilde{f}_2.
\end{equation*}
On $E_1$ the reverser is given by
\begin{equation*}
S:(A,B,\bar{A},\bar{B},q_1,p_1,q_2,p_2)\mapsto (\bar{A},-\bar{B},A,-B,-q_1,p_1,-q_2,p_2).
\end{equation*}
The reduced Hamiltonian $H^\mu$ is independent of $q_1$, $q_2$ since these are cyclic and $p_1$, $p_2$ are therefore preserved. Applying the usual normal form theory for Hamiltonian systems (see \cite{EL}) we may, for every $n_0\geq 2$, write
\begin{equation*}
H^\mu(A,B,\bar{A},\bar{B},0,0)=\i s(A\bar{B}-\bar{A}B)+\abs{B}^2+H_{NF}^0(\abs{A}^2,\i (A\bar{B}-\bar{A}B),\mu)+\mathcal{O}(\abs{(A,B)}^2\abs{(A,B,\mu)}^{n_0}),
\end{equation*}
where $H_{NF}^0$ is real polynomial of degree $n_0+1$ such that $H_{NF}^0(\abs{A}^2,\i(A\bar{B}-\bar{A}B),\mu)=\mathcal{O}(\abs{(A,B)}^2\abs{(A,B,\mu)})$. After a canonical change of variables (see Theorem \ref{normalform})
\begin{align*}
H^\mu(A,B,\bar{A},\bar{B},p_1,p_2)&=\i s(A\bar{B}-\bar{A}B)+\abs{B}^2+\frac{\rho\cos^2(\theta_1)(\alpha-\rho\cos^2(\theta_1))}{2h\alpha\left(\alpha-\cos^2(\theta_1)\left(\rho+\frac{1}{h}\right)\right)}p_1^2\\
&\quad+\frac{\rho\cos^4(\theta_1)}{h\alpha\left(\alpha-\cos^2(\theta_1)\left(\rho+\frac{1}{h}\right)\right)}p_1p_2+\frac{\cos^2(\theta_1)\left(\alpha-\frac{\cos^2(\theta_1)}{h}\right)}{2\alpha\left(\alpha-\cos^2(\theta_1)\left(\rho+\frac{1}{h}\right)\right)}p_2^2\\
&\quad+H_{nl}^\mu(A,B,\bar{A},\bar{B},p_1,p_2),
\end{align*}
where
\begin{align*}
H_{nl}^\mu(A,B,\bar{A},\bar{B},p_1,p_2)&=H_{NF}(\abs{A}^2,\i(A\bar{B}-\bar{A}B),p_1,p_2,\mu)+H_r(A,B,\bar{A},\bar{B},p_1,p_2,\mu)\\
&\quad +\mathcal{O}(\abs{(A,B,p_1,p_2)}^2\abs{(A,B,p_1,p_2,\mu)}^{n_0}),
\end{align*}
with
\begin{align*}
H_{NF}(\abs{A}^2,\i(A\bar{B}-\bar{A}B),0,0,\mu)&=H_{NF}^0(\abs{A}^2,\i(A\bar{B}-\bar{A}B),\mu),\\
H_{NF}(\abs{A}^2,\i(A\bar{B}-\bar{A}B),p_1,p_2,\mu)&=\mathcal{O}(\abs{(A,B)}^2\abs{(A,B,p_1,p_2,\mu)}),
\end{align*}
and 
\begin{equation*}
H_r(A,B,\bar{A},\bar{B},p_1,p_2,\mu)=\mathcal{O}(\abs{(A,B,p_1,p_2)}\abs{(p_1,p_2)}\abs{(p_1,p_2,\mu)})
\end{equation*}
Let $\mu^nH_i^n(A,B,\bar{A},\bar{B},p_1,p_2)$ denote the part of $H^\mu(A,B,\bar{A},\bar{B},p_1,p_2)$ that is homogeneous of order $i$ in $(A,B,\bar{A},\bar{B},p_1,p_2)$ and of order $n$ in $\mu$. We then have that
\begin{align}
H_2^1(A,B,\bar{A},\bar{B},p_1,p_2)&=c_1^1p_1^2+c_2^1\abs{A}^2+c_3^1\i(A\bar{B}-\bar{A}B)+c_4^1p_1A+\bar{c}_4^1p_1\bar{A}+c_5^1p_1B+\bar{c}_5^1p_1\bar{B}\nonumber\\
&\quad +c_7^1p_2^2+c_8^1p_2A+\bar{c}_8^1p_2\bar{A}+c_9^1p_2B+\bar{c}_9^1p_2\bar{B}+c_{10}^1p_1p_2,\label{o21}\\
H_3^0(A,B,\bar{A},\bar{B},p_1,p_2)&=c_1^0p_1^3+c_2^0p_1\abs{A}^2+c_3^0p_1\i(A\bar{B}-\bar{A}B)+c_4^0p_1^2A+\bar{c}_4^0p_1^2\bar{A}+c_5^0p_1^2B+\bar{c}_5^0p_1^2\bar{B}\nonumber\\
&\quad +c_6^0p_2^3+c_7^0p_2\abs{A}^2+c_8^0p_2\i(A\bar{B}-\bar{A}B)+c_9^0p_2^2A+\bar{c}_9^0p_2^2\bar{A}+c_{10}^0p_2^2B+\bar{c}_{10}^0p_2^2\bar{B}\nonumber\\
&\quad+c_{11}^0p_1^2p_2+c_{12}^0p_1p_2^2,\label{o30}\\
H_4^0(A,B,\bar{A},\bar{B},p_1,p_2)&=d_1^0\abs{A}^4+d_2^0\i (A\bar{B}-\bar{A}B)\abs{A}^2-d_3^0(A\bar{B}-\bar{A}B)^2+\mathcal{O}(\abs{(p_1,p_2)}^2\abs{(A,B)}^2)\nonumber\\
&\quad+\mathcal{O}(\abs{(p_1,p_2)}^3\abs{(A,B,p_1,p_2)})\quad\label{o40}
\end{align}
We are interested in the lower order reduced Hamilton's equations;
\begin{align}
A_x&=\frac{\partial H^\mu}{\partial \bar{B}}(A,B,\bar{A},\bar{B},0,0)\label{hameqA},\\
B_x&=-\frac{\partial H^\mu}{\partial \bar{A}}(A,B,\bar{A},\bar{B},0,0)\label{hameqB}.
\end{align}
Using the expansions \eqref{o21}--\eqref{o40}, we find that \eqref{hameqA}--\eqref{hameqB} are given by
\begin{align}
A_x&=\i sA+B+\i c_3^1\mu A+d_2^0\i A\abs{A}^2-2d_3^0A(A\bar{B}-\bar{A}B)+\mathcal{O}(\abs{(A,B)}\abs{(A,B,\mu)}^3)\label{exphameqA},\\
B_x&=\i sB+\i c_3^1\mu B-c_2^1\mu A-2d_1^0A\abs{A}^2-\i d_2^0A^2\bar{B}+2d_2^0B\abs{A}^2-2d_3^0B(A\bar{B}-\bar{A}B)\nonumber\\
&\quad +\mathcal{O}(\abs{(A,B)}\abs{(A,B,\mu)}^3).\label{exphameqB}
\end{align}
We have the following general result regarding reversible systems of the type \eqref{exphameqA}--\eqref{exphameqB}.
\begin{theorem}\label{darkbright}
Suppose that $c_2^1<0$.
\begin{enumerate}
\item \cite{IP} $d_1^0>0$: For each sufficiently small, positive value of $\mu$ the system \eqref{exphameqA}--\eqref{exphameqB} has two distinct symmetric homoclinic solutions.
\item \cite{BG} $d_1^0>0$: For each sufficiently small, positive value of $\mu$ the system \eqref{exphameqA}--\eqref{exphameqB} has two one-parameter families of geometrically distinct homoclinic solutions which generically resemble multiple copies of one of the homoclinic solutions in $1$.
\item \cite{IP} $d_1^0<0$: For each sufficiently small, negative value of $\mu$ the system \eqref{exphameqA}--\eqref{exphameqB} has a one-parameter family of pairs of reversible homoclinic orbits to periodic orbits.
\end{enumerate}
The homoclinic solutions in $1$ and $2$ correspond to travelling waves of amplitude $\mathcal{O}(\mu^\frac{1}{2})$ which have a bright solitary wave profile in the $x$ direction and are $2\pi/(v_0+\mu)$-periodic in $z$. The solutions found in $3$ correspond to travelling waves which have a dark solitary wave profile in the $x$ direction and are $2\pi/(v_0+\mu)$-periodic in $z$. See Figure \ref{darkbrightwaveprofiles} for sketches of the solitary wave profiles in the $x$-direction.
\end{theorem}
In our case 
\begin{equation*}
c_2^1=\frac{2\tilde{\gamma}_1(s\sin(\theta_1)+\nu_0\sin(\theta_2))\sin(\theta_1-\theta_2)}{(s\cos(\theta_1)+\nu_0\cos(\theta_2))(s+\nu_0\cos(\theta_1-\theta_2))\tau_1}\left(\frac{\rho}{\tanh(\tilde{\gamma}_1)}+\frac{1}{\tanh(h\tilde{\gamma}_1)}\right),
\end{equation*}
\begin{align*}
d_1^0&=\frac{-\tilde{\gamma}_1^4}{(s\cos(\theta_1)+\nu_0\cos(\theta_2))^4\tau_1^2}\Bigg\{(s\cos(\theta_1)+\nu_0\cos(\theta_2))^2\bigg[\rho\bigg(-\frac{4\tilde{\gamma}_1}{\tanh(2\tilde{\gamma}_1)\tanh^2(\tilde{\gamma}_1)}+\frac{6\tilde{\gamma}_1}{\tanh(\tilde{\gamma}_1)}\bigg)\\
&\quad-\frac{4\tilde{\gamma}_1}{\tanh(2h\tilde{\gamma}_1)\tanh^2(h\tilde{\gamma}_1)}+\frac{6\tilde{\gamma}_1}{\tanh(h\tilde{\gamma}_1)}\bigg]-\left[\frac{4(s\cos(\theta_1)+\nu_0\cos(\theta_2))^2(s+\nu_0\cos(\theta_1-\theta_2))^2}{\tilde{\gamma}_1^2}\right]\\
&\quad\times\bigg[\frac{\rho}{\tanh^2(\tilde{\gamma}_1)}+\frac{1}{h\tanh^2(h\tilde{\gamma}_1)}\bigg]-\frac{3\tilde{\gamma}_1^4\beta}{2}-\left[\frac{(s\cos(\theta_1)+\nu_0\cos(\theta_2))^4}{2}\right]\\
&\quad\times\bigg[\frac{4}{\tanh(2h\tilde{\gamma}_1)\tanh(h\tilde{\gamma}_1)}+\frac{1}{\sinh^2(h\tilde{\gamma}_1)}-2 -\rho\bigg(\frac{4}{\tanh(2\tilde{\gamma}_1)\tanh(\tilde{\gamma}_1)}+\frac{1}{\sinh^2(\tilde{\gamma}_1)}-2\bigg)\bigg]^2\\
&\quad\times\bigg[\alpha+4\beta\tilde{\gamma}_1^2-\frac{2(s\cos(\theta_1)+\nu_0\cos(\theta_2))^2}{\tilde{\gamma}_1}\bigg(\frac{\rho}{\tanh(2\tilde{\gamma}_1)}+\frac{1}{\tanh(2h\tilde{\gamma}_1)}\bigg)\bigg]^{-1}\\
&\quad-(s\cos(\theta_1)+\nu_0\cos(\theta_2))^2\bigg[\rho\bigg(\frac{s\cos(\theta_1)+\nu_0\cos(\theta_2)}{\sinh^2(\tilde{\gamma}_1)}+\frac{2\cos(\theta_1)(s+\nu_0\cos(\theta_1-\theta_2))}{\tilde{\gamma}_1\tanh(\tilde{\gamma}_1)}\bigg)\\
&\quad-\bigg(\frac{s\cos(\theta_1)+\nu_0\cos(\theta_2)}{\sinh^2(h\tilde{\gamma}_1)}+\frac{2\cos(\theta_1)(s+\nu_0\cos(\theta_1-\theta_2))}{h\tilde{\gamma}_1\tanh(h\tilde{\gamma}_1)}\bigg)\bigg]^2\bigg[\alpha-\cos^2(\theta_1)\bigg(\rho+\frac{1}{h}\bigg)\bigg]^{-1}\Bigg\}.
\end{align*}
It is possible to choose the parameters such that $c_2^1<0,\ d_1^0>0$ or $c_2^1<0,\ d_1^0<0$. This is expected since the coefficients appearing in the Hamiltonian-Hopf bifurcation in $2$-dimensional setting, see \cite{nilsson2016}, satisfies the same property. This bifurcation was investigated for surface waves in \cite{Groves2003} and the results are described in Theorem 6 of that paper. In particular, no dark solitary waves are found and this is due to the fact that the coefficient corresponding to $d_1^0$ is strictly positive when considering surface waves.
{\bf Acknowledgments.} The author was supported by Grant No. 621-2012-3753 from the Swedish Research Council.

The author would also like to thank Erik Wahlén, Mark Groves and Mariana Haragus for their help and advice when writing this article.

\normalsize
\bibliographystyle{plain}

\begin{thebibliography}{10}

\bibitem{Akers2017}
B.~F. Akers and J.~A. Reeger.
\newblock {Three-dimensional overturned traveling water waves}.
\newblock {\em Wave Motion}, 68:210--217, 2017.

\bibitem{Ambrosetti}
A.~Ambrosetti and G.~Prodi.
\newblock {\em {A primer of nonlinear analysis}}.
\newblock Cambridge University Press, Cambridge, 1995.

\bibitem{Bagri2014}
G.~S. Bagri and M.~D. Groves.
\newblock {A Spatial Dynamics Theory for Doubly Periodic Travelling
  Gravity-Capillary Surface Waves on Water of Infinite Depth}.
\newblock {\em J. Dyn. Differ. Equations}, pages 1--28, 2014.

\bibitem{Berti2007}
M.~Berti.
\newblock {\em {Nonlinear oscillations of Hamiltonian PDEs}}.
\newblock Birkh{\"{a}}user Boston, Inc., Boston, MA, 2007.

\bibitem{BG}
B.~Buffoni and M.~D. Groves.
\newblock {A Multiplicity Result for Solitary Gravity-Capillary Waves in Deep
  Water via Critical-Point Theory}.
\newblock {\em Arch. Ration. Mech. Anal.}, 146(3):183--220, 1999.

\bibitem{BGT}
B.~Buffoni, M.~D. Groves, and J.~F. Toland.
\newblock {A plethora of solitary gravity-capillary water waves with nearly
  critical Bond and Froude numbers}.
\newblock {\em Philos. Trans. R. Soc. London. Ser. A. Math. Phys. Sci. Eng.},
  354:575--607, 1996.

\bibitem{Craig2000}
W.~Craig and D.~P. Nicholls.
\newblock {Traveling two and three dimensional capillary gravity water waves}.
\newblock {\em SIAM J. Math. Anal.}, 32(2):323--359, 2000.

\bibitem{EL}
C.~Elphick.
\newblock {Global Aspects of Hamiltonian Normal Forms}.
\newblock {\em Phys. Lett. A}, 127(8,9):418--424, 1988.

\bibitem{gols85}
M.~Golubitsky and D.~G. Schaeffer.
\newblock {\em {Singularities and groups in bifurcation theory. Vol. I}}.
\newblock Springer-Verlag, New York, 1985.

\bibitem{Groves2017}
M.~Groves and D.~Nilsson.
\newblock {Spatial dynamics methods for solitary waves on a ferrofluid jet}.
\newblock {\em arXiv:1706.00453}, 2017.

\bibitem{Groves2001}
M.~D. Groves.
\newblock {An existence theory for three-dimensional periodic travelling
  gravity-capillary water waves with bounded transverse profiles}.
\newblock {\em Phys. D Nonlinear Phenom.}, 152/153:395--415, 2001.

\bibitem{Groves2003}
M.~D. Groves and M.~Haragus.
\newblock {A bifurcation theory for three-dimensional oblique travelling
  gravity-capillary water waves}.
\newblock {\em J. Nonlinear Sci.}, 13(4):397--447, 2003.

\bibitem{Groves2007}
M.~D. Groves and A.~Mielke.
\newblock {A spatial dynamics approach to three-dimensional gravity-capillary
  steady water waves}.
\newblock {\em Proc. R. Soc. Edinburgh Sect. A Math.}, 131:83--136, 2001.

\bibitem{HK2001}
M.~Haragus-Courcelle and K.~Kirchg{\"{a}}ssner.
\newblock {Three-dimensional steady capillary-gravity waves}.
\newblock In {\em Ergod. theory, Anal. Effic. Simul. Dyn. Syst.}, pages
  363----397. Springer, Berlin, 2001.

\bibitem{HaragusCourcelle2000}
M.~Haragus-Courcelle and R.L. Pego.
\newblock {Spatial wave dynamics of steady oblique wave interactions}.
\newblock {\em Phys. D Nonlinear Phenom.}, 145(3-4):207--232, 2000.

\bibitem{IP}
G.~Iooss and M.~P{\'{e}}rou{\`{e}}me.
\newblock {Perturbed homoclinic solutions in reversible 1:1 resonance vector
  fields}.
\newblock {\em J. Differ. Equ.}, 102(1):62--88, 1993.

\bibitem{IP3d}
G.~Iooss and P.~Plotnikov.
\newblock {Small divisor problem in the theory of three-dimensional water
  gravity waves}.
\newblock {\em Mem. Amer. Math. Soc.}, 200:viii+128, 2009.

\bibitem{KA}
T.~Kato.
\newblock {\em {Perturbation theory for linear operators}}.
\newblock Springer-Verlag New York, Inc., New York, 1966.

\bibitem{KIM2006}
B.~Kim and T.~R. Akylas.
\newblock {On gravity–capillary lumps. Part 2. Two-dimensional Benjamin
  equation}.
\newblock {\em J. Fluid Mech.}, 557:237, 2006.

\bibitem{K82}
K.~Kirchg{\"{a}}ssner.
\newblock {Wave-Solutions of Reversible Systems and Applications}.
\newblock {\em J. Differ. Equations}, 45:113--127, 1982.

\bibitem{MI}
A.~Mielke.
\newblock {Reduction of quasilinear elliptic equations in cylindrical domains
  with applications}.
\newblock {\em Math. Methods Appl. Sci.}, 10:51--66, 1988.

\bibitem{nilsson2016}
D.~Nilsson.
\newblock {Internal gravity-capillary solitary waves in finite depth}.
\newblock {\em Math. Methods Appl. Sci.}, 40(4):1053--1080, 2017.

\bibitem{Parau2007}
E.~Parau, J-M Vanden-Broeck, and M.~Cooker.
\newblock {Nonlinear three-dimensional interfacial flows with a free surface}.
\newblock {\em J. Fluid Mech.}, 591:481--494, 2007.

\bibitem{Parau2007a}
E.~Parau, J-M. Vanden-Broeck, and M.~Cooker.
\newblock {Three-dimensional gravity and gravity-capillary interfacial flows}.
\newblock {\em Math. Comput. Simul.}, 74:105--112, 2007.

\bibitem{Reeder1981}
J.~Reeder and M.~Shinbrot.
\newblock {Three-dimensional, nonlinear wave interaction in water of constant
  depth}.
\newblock {\em Nonlinear Anal.}, 5(3):303--323, 1981.

\end{thebibliography}

\end{document}